\documentclass[a4paper,twoside]{article}
\usepackage{amssymb}

\usepackage{amsbsy}
\usepackage{amsmath}
\usepackage{theorem}
\usepackage{graphicx}
\usepackage{epsfig}

\evensidemargin\oddsidemargin
\numberwithin{equation}{section}
\newtheorem{theorem}{Theorem}[section]
\newtheorem{corollary}[theorem]{Corollary}

\newtheorem{lemma}[theorem]{Lemma}
\newtheorem{proposition}[theorem]{Proposition}

{\theorembodyfont{\rmfamily}
\newtheorem{definition}[theorem]{Definition}
}
{\theorembodyfont{\rmfamily}
\newtheorem{remark}[theorem]{Remark}
}
{\theorembodyfont{\rmfamily}

}
{\theorembodyfont{\rmfamily}
\newtheorem{example}[theorem]{Example}
}
{\theorembodyfont{\rmfamily}

}
{\theorembodyfont{\rmfamily}
\newtheorem{final remarks and questions}[theorem]{Final remarks and questions}
}
\begin{document}

\title{\textbf{Resolving 3-dimensional toric singularities}}
\author{\thanks{{\scriptsize GSRT-fellow supported by the European Union and the
Greek Ministry of Research and Technology.}\medskip\newline
} \textbf{Dimitrios
I. Dais} \\
\noindent {\scriptsize Mathematics Department, Section of Algebra and
Geometry, University of Ioannina}\\
{\scriptsize GR-45110 Ioannina, Greece. E-mail: ddais@cc.uoi.gr}}
\date{}
\maketitle

\begin{abstract}
\noindent This paper surveys, in the first place, some basic facts from the
classification theory of normal complex singularities, including details for
the low dimensions $2$ and $3$. Next, it describes how the toric
singularities are located within the class of rational singularities, and
recalls their main properties. Finally, it focuses, in particular, on a
toric version of Reid's desingularization strategy in dimension three.
\end{abstract}

\section{Introduction\label{SECTION1}}

\noindent There are certain general \textit{qualitative criteria} available
for the rough classification of singularities of complex varieties. The main
ones arise:\bigskip
\begin{equation*}
\left\{
\begin{array}{l}
\bullet \text{ from the study of the punctual algebraic} \\
\text{ behaviour of these varieties} \\
\text{(w.r.t. local rings associated to singular points) } \\
\text{\lbrack \textit{algebraic classification}]} \\
\, \\
\bullet \text{ from an intrinsic characterization } \\
\text{for the nature of the possible exceptional } \\
\text{loci w.r.t. any desingularization} \\
\text{\lbrack \textit{rational, elliptic, non-elliptic etc.}]} \\
\, \\
\bullet \text{ from the behaviour of ``discrepancies''} \\
\text{(for }\mathbb{Q}\text{-Gorenstein normal complex varieties)} \\
\text{\lbrack \textit{adjunction-theoretic classification}]}
\end{array}
\right.
\end{equation*}
\newpage

\noindent {}$\bullet $ \textbf{Algebraic Classification.} At first we recall
some fundamental definitions from commutative algebra (cf. \cite{Kunz}, \cite
{Matsumura}). Let $R$ be a commutative ring with $1$. The \textit{height} ht$%
\left( \frak{p}\right) $ of a prime ideal $\frak{p}$ of $R$ is the supremum
of the lengths of all prime ideal chains which are contained in $\frak{p}$,
and the \textit{dimension} of $R$ is defined to be
\begin{equation*}
\text{dim}\left( R\right) :=\text{sup}\left\{ \text{ht}\left( \frak{p}%
\right) \left| \frak{p}\text{ prime ideal of }R\right. \right\} .
\end{equation*}
$R$ is \textit{Noetherian} if any ideal of it has a finite system of
generators. $R$ is a \textit{local ring} if it is endowed with a \textit{%
unique} maximal ideal $\frak{m}$. A local ring $R$ is \textit{regular }%
(resp. \textit{normal}) if dim$\left( R\right) =$ dim$\left( \frak{m}/\frak{m%
}^{2}\right) $ (resp. if it is an integral domain and is integrally closed
in its field of fractions). A finite sequence $a_{1},\ldots ,a_{\nu }$ of
elements of a ring $R$ is defined to be a \textit{regular sequence} if $a_{1}
$ is not a zero-divisor in $R$ and for all $i$, $i=2,\ldots ,\nu $, $a_{i}$
is not a zero-divisor of $R/\left\langle a_{1},\ldots ,a_{i-1}\right\rangle $%
. A Noetherian local ring $R$ (with maximal ideal $\frak{m}$) is called
\textit{Cohen-Macaulay} if
\begin{equation*}
\text{depth}\left( R\right) =\text{dim}\left( R\right) ,
\end{equation*}
where the \textit{depth} of $R$ is defined to be the maximum of the lengths
of all regular sequences whose members belong to $\frak{m}$. A
Cohen-Macaulay local ring $R$ is called \textit{Gorenstein }if
\begin{equation*}
\text{Ext}_{R}^{\text{dim}\left( R\right) }\left( R/\frak{m},R\right) \cong
R/\frak{m}.
\end{equation*}
A Noetherian local ring $R$ is said to be a \textit{complete intersection}
if there exists a regular local ring $R^{\prime }$, such that $R\cong
R^{\prime }/\left( f_{1},\ldots ,f_{q}\right) $ for a finite set $\left\{
f_{1},\ldots ,f_{q}\right\} \subset R^{\prime }$ whose cardinality equals $q=
$ dim$\left( R^{\prime }\right) -$ dim$\left( R\right) $. The hierarchy by
inclusion of the above types of Noetherian local rings is known to be
described by the following diagram:
\begin{equation}
\begin{array}{ccc}
\left\{ \text{Noetherian local rings}\right\}  & \supset  & \left\{ \text{%
normal local rings}\right\}  \\
\cup  &  & \cup  \\
\left\{ \text{Cohen-Macaulay local rings}\right\}  &  & \left\{ \text{%
regular local rings}\right\}  \\
\cup  &  & \cap  \\
\left\{ \text{Gorenstein local rings}\right\}  & \supset  & \left\{ \text{%
complete intersections (``c.i.'s'')}\right\}
\end{array}
\label{HIERARCHY}
\end{equation}
An arbitrary Noetherian ring $R$ and its associated affine scheme Spec$%
\left( R\right) $ are called Cohen-Macaulay, Gorenstein, normal or regular,
respectively, iff \ all the localizations $R_{\frak{m}}$ with respect to all
the members $\frak{m}\in $ Max-Spec$\left( R\right) $ of the maximal
spectrum of $R$ are of this type. In particular, if the $R_{\frak{m}}$'s for
all maximal ideals $\frak{m}$ of $R$ are c.i.'s, then one often says that $R$
is a \textit{locally complete intersection }(``l.c.i.'') to distinguish it
from the ``global'' ones. (A \textit{global complete intersection}
(``g.c.i.'') is defined to be a ring $R$ of finite type over a field $%
\mathbf{k}$ (i.e., an affine $\mathbf{k}$-algebra), such that \
\begin{equation*}
R\cong \mathbf{k}\left[ \mathsf{T}_{1}..,\mathsf{T}_{d}\right] \,/\,\left(
\varphi _{1}\left( \mathsf{T}_{1},..,\mathsf{T}_{d}\right) ,..,\varphi
_{q}\left( \mathsf{T}_{1},..,\mathsf{T}_{d}\right) \right)
\end{equation*}
for $q$ polynomials $\varphi _{1},\ldots ,\varphi _{q}$ from $\mathbf{k}%
\left[ \mathsf{T}_{1},..,\mathsf{T}_{d}\right] $ with $q=d-$ dim$\left(
R\right) $). Hence, the above inclusion hierarchy can be generalized for all
Noetherian rings, just by omitting in (\ref{HIERARCHY}) the word ``local''
and by substituting l.c.i.'s for c.i.'s.\medskip

We shall henceforth consider only \textit{complex varieties }$\left( X,%
\mathcal{O}_{X}\right) $, i.e., integral separated schemes of finite type
over $\mathbf{k}=\mathbb{C}$; thus, the punctual algebraic behaviour of $X$
is determined by the stalks $\mathcal{O}_{X,x}$ of its structure sheaf $%
\mathcal{O}_{X}$, and $X$ itself is said to have a given \textit{algebraic}
\textit{property} whenever all $\mathcal{O}_{X,x}$'s have the analogous
property from (\ref{HIERARCHY}) for all $x\in X$. Furthermore, via the
\textsc{gaga}-correspondence (\cite{Serre}, \cite[\S 2]{Gro2}) which
preserves the above quoted algebraic properties, we may work within the
\textit{analytic category} by using the usual contravariant functor
\begin{equation*}
\left( X,x\right) \rightsquigarrow \mathcal{O}_{X,x}^{\text{hol}}
\end{equation*}
between the category of isomorphy classes of germs of $X$ and the
corresponding category of isomorphy classes of analytic local rings at the
marked points $x$. For\textsf{\ }a complex variety $X$ and $x\in X$, we
denote by $\frak{m}_{X,x}$ the maximal ideal of $\mathcal{O}_{X,x}^{\text{hol%
}}$ and by
\begin{align}
\text{Sing}\left( X\right) & =\left\{ x\in X\ \left| \ \right. \mathcal{O}%
_{X,x}^{\text{hol}}\text{ is a non-regular local ring}\right\} \medskip
\label{SING-LOC} \\
& =\left\{ x\in X\ \left| \ \right. \dim \left( \frak{m}_{X,x}/\frak{m}%
_{X,x}^{2}\right) >\dim _{x}\left( X\right) \right\}  \notag
\end{align}
its \textit{singular locus}. By a \textit{desingularization} (or \textit{%
resolution of singularities}) $f:\widehat{X}\rightarrow X$ of a non-smooth $%
X $, we mean a ``full'' or ``overall'' desingularization (if not mentioned),
i.e., Sing$(\widehat{X})=\varnothing $. When we deal with \textit{partial}
desingularizations, we mention it explicitly.\bigskip

\noindent {}$\bullet $ \textbf{Rational and Elliptic Singularities}. We say
that $X$ has (at most) \textit{rational singularities} if there exists a
desingularization $f:Y\rightarrow X$ of $X$, such that
\begin{equation*}
f_{\ast }\mathcal{O}_{Y}=\mathcal{O}_{X}
\end{equation*}
(equivalently, $Y$ is normal), and
\begin{equation*}
R^{i}f_{\ast }\mathcal{O}_{Y}=0,\ \ \ \ \forall i,\ \ \ 1\leq i\leq \text{
dim}_{\mathbb{C}}X-1.
\end{equation*}
(The $i$-th direct image sheaf is defined via
\begin{equation*}
U\longmapsto R^{i}f_{\ast }\mathcal{O}_{Y}\left( U\right) :=H^{i}\left(
f^{-1}\left( U\right) ,\mathcal{O}_{Y}\left| _{f^{-1}\left( U\right)
}\right. \right) \ ).
\end{equation*}
This definition is independent of the particular choice of the
desingularization of $X.$ (Standard example: quotient singularities\footnote{%
The \textit{quotient singularities} are of the form $(\mathbb{C}^{r}/G,[%
\mathbf{0}]),$ where $G$ is a finite subgroup of GL$(r,\mathbb{C})$ (without
pseudoreflections) acting linearly on $\mathbb{C}^{r}$, $p:$ $\mathbb{C}%
^{r}\rightarrow \mathbb{C}^{r}/G=$ Spec$(\mathbb{C}[\mathbf{z}]^{G})$ the
quotient map, and $[\mathbf{0}]=p(\mathbf{0}).$ Note that
\begin{equation*}
\text{Sing}\left( \mathbb{C}^{r}/G\right) =p\left( \left\{ \mathbf{z}\in
\mathbb{C}^{r}\ \left| \ G_{\mathbf{z}}\neq \left\{ \text{Id}\right\}
\right. \right\} \right)
\end{equation*}
(cf. (\ref{SING-LOC})), where $G_{\mathbf{z}}:=\left\{ g\in G\ \left| \
g\cdot \mathbf{z=z}\right. \right\} $ is the isotropy group of $\mathbf{z}%
\in \mathbb{C}^{r}$.} are rational singularities).\medskip

\noindent {}We say that a Gorenstein singularity $x$ of $X$ is an \textit{%
elliptic singularity} if there exists a desingularization $f:Y\rightarrow X$
of $x\in X$, such that
\begin{equation*}
R^{i}f_{\ast }\mathcal{O}_{Y}=0,\ \ \forall i,\ \ \ 1\leq i\leq \text{dim}_{%
\mathbb{C}}X-2,\ \ \text{ }
\end{equation*}
and
\begin{equation*}
R^{\dim _{\mathbb{C}}X-1}f_{\ast }\mathcal{O}_{Y}\cong \mathbb{C}.
\end{equation*}
(The definition is again independent of the particular choice of the
desingularization).\newpage

\bigskip

\noindent {}$\bullet $ \textbf{Adjunction-Theoretic Classification. }If $X$
is a \textit{normal} complex variety, then its Weil divisors can be
described by means of ``divisorial'' sheaves as follows:

\begin{lemma}
\emph{(\cite[1.6]{HART3}). }For a coherent sheaf $\mathcal{F}$ of $\mathcal{O%
}_{X}$-modules the following conditions are equivalent\emph{:}\newline
\emph{(i)} $\mathcal{F}$ is \emph{reflexive} \emph{(}i.e., $\mathcal{F}\cong
\mathcal{F}^{\vee \vee },$ with $\mathcal{F}^{\vee }:=Hom_{\mathcal{O}%
_{X}}\left( \mathcal{F},\mathcal{O}_{X}\right) $ denoting the dual of $%
\mathcal{F}$\emph{) }and has rank one.\newline
\emph{(ii)} If $X^{0}$ is a non-singular open subvariety of $X$ with \emph{%
codim}$_{X}\left( X\mathbb{r}X^{0}\right) \geq 2,$ then $\mathcal{F}\left|
_{X^{0}}\right. $ is invertible and
\begin{equation*}
\mathcal{F}\cong \iota _{\ast }\left( \mathcal{F}\left| _{X^{0}}\right.
\right) \cong \iota _{\ast }\iota ^{\ast }\left( \mathcal{F}\right) ,
\end{equation*}
where $\iota :X^{0}\hookrightarrow X$ denotes the inclusion map.
\end{lemma}

\noindent {}The \textit{divisorial sheaves} are exactly those satisfying one
of the above conditions. Since a divisorial sheaf is torsion free, there is
a non-zero section $\gamma \in H^{0}\left( X,Rat_{X}\otimes _{\mathcal{O}%
_{X}}\mathcal{F}\right) $, with
\begin{equation*}
H^{0}\left( X,Rat_{X}\otimes _{\mathcal{O}_{X}}\mathcal{F}\right) \cong
\mathbb{C}\left( X\right) \cdot \gamma ,
\end{equation*}
and $\mathcal{F}$ can be considered as a subsheaf of the constant sheaf $%
Rat_{X}$ of rational functions of $X,$ i.e., as a special \textit{fractional
ideal sheaf.}

\begin{proposition}
\emph{(\cite[App. of \S 1]{REID1}) }The correspondence
\begin{equation*}
Cl\left( X\right) \ni \{D\}\overset{\delta }{\longmapsto }\{\mathcal{O}%
_{X}\left( D\right) \}\in \left\{
\begin{array}{c}
\text{divisorial coherent } \\
\text{subsheaves of }Rat_{X}
\end{array}
\right\} \ /\ H^{0}\left( X,\mathcal{O}_{X}^{\ast }\right)
\end{equation*}
with $\mathcal{O}_{X}\left( D\right) $ defined by sending every non-empty
open set $U$ of $X$ onto
\begin{equation*}
U\longmapsto \mathcal{O}_{X}\left( D\right) \left( U\right) :=\left\{
\varphi \in \mathbb{C}\left( X\right) ^{\ast }\ \left| \ \left( \text{\emph{%
div}}\left( \varphi \right) +D\right) \left| _{U}\right. \geq 0\right.
\right\} \cup \{0\},
\end{equation*}
is a bijection, and induces a $\mathbb{Z}$-module isomorphism. In fact, to
avoid torsion, one defines this $\mathbb{Z}$-module structure by setting
\begin{equation*}
\delta \left( D_{1}+D_{2}\right) :=(\mathcal{O}_{X}\left( D_{1}\right)
\otimes \mathcal{O}_{X}\left( D_{2}\right) )^{\vee \vee }\text{\emph{and}\ }%
\delta \left( \kappa D\right) :=\mathcal{O}_{X}\left( D\right) ^{\left[
\kappa \right] }=\mathcal{O}_{X}\left( \kappa D\right) ^{\vee \vee },
\end{equation*}
for any Weil divisors $D,D_{1},D_{2}$ and $\kappa \in \mathbb{Z}.$
\end{proposition}

\noindent {}Let now $\Omega _{\text{Reg}(X)/\mathbb{C}}$ be the sheaf of
regular $1$-forms, or K\"{a}hler differentials, on
\begin{equation*}
\text{Reg}\left( X\right) =X\mathbb{r}\text{Sing}\left( X\right) \overset{%
\iota }{\hookrightarrow }X,
\end{equation*}
(cf. \cite[\S 5.3]{IITAKA}) and for $i\geq 1,$ let us set
\begin{equation*}
\Omega _{\text{Reg}(X)/\mathbb{C}}^{i}:=\bigwedge\limits^{i}\ \Omega _{\text{%
Reg}(X)/\mathbb{C}}\ .
\end{equation*}
The unique (up to rational equivalence) Weil divisor $K_{X},$ which maps
under $\delta $ to the \textit{canonical divisorial sheaf}
\begin{equation*}
\omega _{X}:=\iota _{\ast }\left( \Omega _{\text{Reg}(X)/\mathbb{C}}^{\dim _{%
\mathbb{C}}\left( X\right) }\right) ,
\end{equation*}
is called the \textit{canonical divisor} of $X.$ Another equivalent
interpretation of $\omega _{X},$ when $X$ is Cohen-Macaulay, can be given by
means of the Duality Theory (see \cite{HART1}, \cite{Gro1}). If $\mathbb{D}%
_{c}^{+}\left( \mathcal{O}_{X}\right) $ denotes the derived category of
below bounded complexes whose cohomology sheaves are coherent, then there
exists a \textit{dualizing complex}\footnote{%
There is a canonical morphism $\omega _{X}\left[ \dim _{\mathbb{C}}\left(
X\right) \right] \rightarrow \omega _{X}^{\bullet }$ which is a
quasi-isomorphism iff $X$ is Cohen-Macaulay.} $\omega _{X}^{\bullet }$ $\in
\mathbb{D}_{c}^{+}\left( \mathcal{O}_{X}\right) $ over $X.$ If $X$ is
Cohen-Macaulay, then the $i$-th cohomology sheaf $\mathcal{H}^{i}\left(
\omega _{X}^{\bullet }\right) $ vanishes for all $i\in \mathbb{Zr}\{-\dim _{%
\mathbb{C}}\left( X\right) \},$ and $\omega _{X}\cong \mathcal{H}^{-\dim _{%
\mathbb{C}}\left( X\right) }\left( \omega _{X}^{\bullet }\right) .$ This
leads to the following:

\begin{proposition}
\label{PROP-GOR}A normal complex variety $X$ is Gorenstein if and only if it
is Cohen-Macaulay and $\omega _{X}$ is invertible.
\end{proposition}

\noindent {}\textit{Proof.} If $X$ is Gorenstein, then $\mathcal{O}_{X,x}$
satisfies the equivalent conditions of \cite[Thm. 18.1]{Matsumura}, for all $%
x\in X.$ This means that $\mathcal{O}_{X,x}$ (as Noetherian local ring) is a
dualizing complex for itself (cf. \cite[Ch. V, Thm. 9.1, p. 293]{HART1}).
Since dualizing complexes are unique up to tensoring with an invertible
sheaf, say $L$, over $X,$ shifted by an integer $n$ (cf. \cite[Ch. V, Cor.
2.3, p. 259]{HART1}), we shall have $\omega _{X}^{\bullet }\cong $ $\mathcal{%
O}_{X,x}^{\bullet }$ $\otimes L\left[ n\right] $. Hence, $\omega _{X}$
itself will be also invertible. The converse follows from the isomorphisms $%
\omega _{X}\cong \mathcal{H}^{-\dim _{\mathbb{C}}\left( X\right) }\left(
\omega _{X}^{\bullet }\right) $ and $\omega _{X,x}\cong \mathcal{O}_{X,x}$,
for all $x\in X.$ (Alternatively, one may use the fact that $x\in X$ is
Gorenstein iff $\mathcal{O}_{X,x}$ is Cohen-Macaulay and $H_{\frak{m}%
_{X,x}}^{\dim _{\mathbb{C}}\left( X\right) }(\mathcal{O}_{X,x})$ is a
dualizing module for it, cf. \cite[Prop. 4.14, p. 65]{Gro1}. The classical
duality \cite[Thm. 6.3, p. 85]{Gro1}, \cite[Ch. V, Cor. 6.5, p. 280]{HART1},
combined with the above uniqueness argument, gives again the required
equivalence).\hfill $\square $

\begin{theorem}
\label{EF}\emph{(Kempf \cite[p. 50]{KKMS}, Elkik \cite{ELKIK1}, \cite{EKLIK2}%
, Bingener-Storch \cite{B-S}). }Let $X$ a normal complex variety of
dimension $\geq 2.$ Then
\begin{equation*}
\left(
\begin{array}{c}
X\text{ \thinspace has at most } \\
\text{rational singularities}
\end{array}
\right) \Longleftrightarrow \left(
\begin{array}{c}
X\text{ \thinspace is Cohen-Macaulay } \\
\text{and }\omega _{X}\cong f_{\ast }\omega _{Y}
\end{array}
\right) ,
\end{equation*}
where $f:Y\longrightarrow X$ is any desingularization of $X.$
\end{theorem}

\noindent {}(Note that, if $E=f^{-1}\left( \text{Sing}\left( X\right)
\right) $ and $\iota :$ Reg$\left( X\right) \hookrightarrow X,\ j:Y\mathbb{r}%
E\hookrightarrow Y$ are the natural inclusions, then by the commutative
diagram
\begin{equation*}
\begin{array}{ccc}
Y\mathbb{r}E & \overset{j}{\hookrightarrow } & Y \\
\downarrow \cong &  & _{f}\downarrow \\
\text{Reg}\left( X\right) & \overset{\iota }{\hookrightarrow } & X
\end{array}
\end{equation*}
we have in general $f_{\ast }\omega _{Y}\hookrightarrow f_{\ast }j_{\ast
}\left( \omega _{Y}\left| _{Y\mathbb{r}E}\right. \right) =\iota _{\ast
}f_{\ast }\left( \omega _{Y}\left| _{Y\mathbb{r}E}\right. \right) =\iota
_{\ast }\left( \omega _{\text{Reg}\left( X\right) }\right) \cong \omega
_{X}. $ In fact, $f_{\ast }\omega _{Y}$ does not depend on the particular
choice of the desingularization.)\medskip

\noindent {}\textit{Sketch of proof}: Let $G_{X}:=$ Co$\ker \left( f_{\ast
}\omega _{Y}\hookrightarrow \omega _{X}\right) ,$ and
\begin{equation*}
\mathcal{L}_{X}^{\bullet }=\{\mathcal{L}_{X}^{i}\}\in \mathbb{D}%
_{c}^{+}\left( \mathcal{O}_{X}\right) ,\ \ \ \ \mathcal{M}_{X}^{\bullet }=\{%
\mathcal{M}_{X}^{i}\}\in \mathbb{D}_{c}^{+}\left( \mathcal{O}_{X}\right) ,
\end{equation*}
the map cones of the canonical homomorphisms
\begin{equation*}
\mathcal{O}_{X}\longrightarrow R^{\bullet }f_{\ast }\mathcal{O}_{Y}\text{ \
\ \ and\ \ \ }f_{\ast }\left( \omega _{Y}[\dim _{\mathbb{C}}\left( X\right)
]\right) \longrightarrow \omega _{X}^{\bullet },
\end{equation*}
respectively. By Grauert-Riemenschneider vanishing theorem, $R^{i}f_{\ast
}\omega _{Y}=0$ for all $i\in \mathbb{Z}_{\geq 1},$ which means that the
canonical morphism
\begin{equation*}
f_{\ast }\left( \omega _{Y}[\dim _{\mathbb{C}}\left( X\right) ]\right)
\longrightarrow \mathbf{R}f_{\ast }\left( \omega _{Y}[\dim _{\mathbb{C}%
}\left( X\right) ]\right)
\end{equation*}
\ is an isomorphism. Moreover, $\mathcal{H}^{i}\left( \mathcal{L}%
_{X}^{\bullet }\right) =\mathcal{L}_{X}^{i}$, for all $\ i\in \mathbb{Z},$%
and
\begin{equation*}
\mathcal{H}^{i}\left( \mathcal{M}_{X}^{\bullet }\right) =\left\{
\begin{array}{ll}
\mathcal{H}^{i}\left( \omega _{X}^{\bullet }\right) , & \text{for }i\in
\lbrack -\dim _{\mathbb{C}}\left( X\right) +1,-t\left( X\right) ], \\
G_{X} & \text{for }i=-\dim _{\mathbb{C}}\left( X\right) , \\
0 & \text{otherwise,}
\end{array}
\smallskip \right.
\end{equation*}
where $t(X)=\inf \{\left. \text{depth}\left( \mathcal{O}_{X,x}\right)
\,\right| \,x\in X\}.$ By the Duality Theorem for proper maps (cf. \cite[Ch.
III, Thm. 11.1, p. 210]{HART1}) we obtain a canonical isomorphism
\begin{equation*}
f_{\ast }\left( \omega _{Y}[\dim _{\mathbb{C}}\left( X\right) ]\right) \cong
\mathbf{R}f_{\ast }\left( \omega _{Y}[\dim _{\mathbb{C}}\left( X\right)
]\right) \text{ }\cong \mathbf{R}\text{\textit{Hom}}_{\mathcal{O}_{X}}\left(
R^{\bullet }f_{\ast }\mathcal{O}_{Y},\omega _{X}^{\bullet }\right)
\end{equation*}
By dualization this reads as
\begin{equation*}
\mathbf{R}\text{\textit{Hom}}_{\mathcal{O}_{X}}\left( \text{ }f_{\ast
}\left( \omega _{Y}[\dim _{\mathbb{C}}\left( X\right) ]\right) ,\omega
_{X}^{\bullet }\right) \cong R^{\bullet }f_{\ast }\mathcal{O}_{Y}.
\end{equation*}
From this isomorphism we deduce that
\begin{equation*}
\mathcal{H}^{i}\left( \mathcal{M}_{X}^{\bullet }\right) =Ext_{\mathcal{O}%
_{X}}^{i+1}\left( \mathcal{M}_{X}^{\bullet },\omega _{X}^{\bullet }\right)
\text{ \ \ and \ \ }\mathcal{H}^{i}\left( \mathcal{L}_{X}^{\bullet }\right)
=Ext_{\mathcal{O}_{X}}^{i+1}\left( \mathcal{L}_{X}^{\bullet },\omega
_{X}^{\bullet }\right)
\end{equation*}
Hence, the assertion follows from the equivalence: $\mathcal{L}_{X}^{\bullet
}=0\Longleftrightarrow \mathcal{M}_{X}^{\bullet }=0.$\hfill$\square $

\begin{remark}
For another approach, see Flenner \cite[Satz 1.3]{FLENER}. For a proof which
avoids Duality Theory, cf. \cite[Cor. 11.9, p. 281]{Kollar} or \cite[Lemma
5.12, p. 156]{Ko-Mo}.
\end{remark}

\begin{definition}
A normal complex variety $X$ is called $\mathbb{Q}$-\textit{Gorenstein} if
\begin{equation*}
\omega _{X}=\mathcal{O}_{X}\left( K_{X}\right)
\end{equation*}
with $K_{X}$ a $\mathbb{Q}$-Cartier divisor. (The smallest positive integer $%
\ell ,$ for which $\ell K_{X}$ is Cartier, is called \textit{the index} of $%
X $.) Let $X$ be a singular $\mathbb{Q}$-Gorenstein complex variety of
dimension $\geq 2.$ Take a desingularization $f:Y\rightarrow X$ of $X$, such
that the exceptional locus of $f$ is a divisor $\bigcup\limits_{i}D_{i}$
with only simple normal crossings, and define \textit{the discrepancy}%
\footnote{%
We may formally define the\textit{\ pull-back} $f^{\ast }\left( K_{X}\right)
$ as the $\mathbb{Q}$-Cartier divisor $\frac{1}{\ell }f^{\ast }\left( \ell
K_{X}\right) $, where $\ell $ is the index of $X.$}
\begin{equation*}
K_{Y}-f^{\ast }\left( K_{X}\right) =\sum_{i}a_{i}D_{i}\ .
\end{equation*}
We say that $X$ has \textit{terminal} (resp. \textit{canonical}, resp.
\textit{log-terminal, }resp.\textit{\ log-canonical}) \textit{singularities}
if all $a_{i}$'s are $>0$ (resp. $\geq 0\ /\ >-1\ /\ \geq -1$). This
definition is independent of the particular choice of the desingularization.
\end{definition}

\begin{remark}
(i) If all $a_{i}$'s are $=0,$ then $f:Y\rightarrow X$ is called a \textit{%
crepant} desingularization of $X$. In fact, the number of crepant divisors $%
\#\left\{ i\ \left| \ a_{i}=0\right. \right\} $ remains invariant w.r.t. all
$f$'s as long as $X$ has at most canonical singularities.\smallskip \newline
(ii) Terminal singularities constitute the smallest class of singularities
to run the MMP (= minimal model program) for smooth varieties. The canonical
singularities are precisely the singularities which appear on the canonical
models of varieties of general type. Finally, the log-singularities are
those singularities for which the \textit{discrepancy function} (assigned to
$\mathbb{Q}$-Gorenstein complex varieties $X$) still makes sense\footnote{%
Cf. \cite[6.3, p. 39]{CKM}, \cite[Prop. 1.9, p. 14]{KOETAL} and \cite[p. 57]
{Ko-Mo}.}. For details about the general MMP, see \cite{KMM}, \cite[6.3, p.
39]{CKM}, \cite{KOETAL} and \cite{Ko-Mo}.
\end{remark}

\begin{theorem}
\label{LOG-RAT}Log-terminal singularities are rational.
\end{theorem}

\noindent {}\textit{Proof}: This follows from \cite[Thm. 1-3-6, p. 311]{KMM}%
, \cite[Cor. 11.14, p. 283]{Kollar} or \cite[Thm. 5.22, p. 161]{Ko-Mo}.\hfill%
$\square $

\begin{corollary}
A singularity $x\in X$ is canonical of index $1$ if and only if it is
rational and Gorenstein.
\end{corollary}

\noindent {}\textit{Proof}: ``$\Rightarrow $'' is obvious by Thm. \ref
{LOG-RAT} and Proposition \ref{PROP-GOR}. (The rationality of canonical
singularities was first shown by Elkik \cite{EKLIK2}). ``$\Leftarrow $''
follows from the fact that $\omega _{X}$ is locally free and from $\omega
_{X}\cong f_{\ast }\omega _{Y}$ (via the other direction of Thm. \ref
{LOG-RAT}), as this is equivalent for $x\in X$ to be canonical of index $1$
(cf. \cite[(1.1), p. 276]{REID1}).\hfill$\square $

\begin{definition}
Let $(\mathcal{O}_{X,x},\frak{m}_{X,x})$ the local ring of a point $x$ of a
normal quasiprojective complex variety $X$ and $V_{x}\subset \mathcal{O}%
_{X,x}$ a finite-dimensional $\mathbb{C}$-vector space generating $\frak{m}%
_{X,x}.$ A \textit{general hyperplane section through} $x$ is a $\mathbb{C}$%
-algebraic subscheme $\mathbb{H}\subset U_{x}$ determined in a suitable
Zariski-open neighbourhood $U_{x}$ of $x$ by the ideal sheaf $\mathcal{O}%
_{X}\cdot v,$ where $v\in V_{x}$ is sufficiently general. (\textit{%
Sufficiently general} means that $v$ can be regarded as a $\mathbb{C}$-point
of a whole Zariski-open dense subset of $V_{x}.)$
\end{definition}

\begin{theorem}
\emph{(M. Reid, \cite[2.6]{REID1}, \cite[3.10]{REID4}, \cite[5.30-5.31, p.
164]{Ko-Mo}).} Let $X$ be a normal quasiprojective complex variety of
dimension $r\geq 3$ and $x\in \,$\emph{Sing}$\left( X\right) .$ If $\left(
X,x\right) $ is a rational Gorenstein singularity, then, for a general
hyperplane section $\mathbb{H}$ through $x,\ \left( \mathbb{H},x\right) $ is
either a rational or an elliptic $(r-1)$-dimensional singularity.
\end{theorem}

\section{Basic facts about two- and three-dimensional \newline
normal singularities\label{SECTION2}}

\noindent {}In dimension 2, the definition of rational and elliptic
singularities fits quite well our intuition of what ``rational'' and
``elliptic'' ought to be. ``Terminal'' points are the smooth ones and the
canonical singularities turn out to be the traditional RDP's (see below
Theorem \ref{RDP}). Moreover, terminal and canonical points have always
index $1.$ On the other hand, the existence of a unique \textit{minimal}%
\footnote{%
A desingularization $f:X^{\prime }\rightarrow X$ of a normal surface $X$ is
minimal if Exc$\left( f\right) $ does not contain any curve with
self-intersection number $-1$ or, equivalently, if for an arbitrary
desingularization $g:X^{\prime \prime }\rightarrow X$ of $X,$ there exists a
unique morphism $h:X^{\prime \prime }\rightarrow X^{\prime }$ with $g=f\circ
h.$} (and \textit{good minimal}\footnote{%
A desingularization of a normal surface is good if (i) the irreducible
components of the exceptional locus are smooth curves, and (ii) the support
of the inverse image of each singular point is a divisor with simple normal
crossings. For the proof of the uniqueness (up to a biregular isomorphism)
of both minimal and good minimal desingularizations, see Brieskorn
\cite[Lemma 1.6]{BRIESKORN1} and Laufer \cite[Thm. 5.12]{Laufer1}.})
desingularization makes the study of normal surface singularities easier
that in higher dimensions.

\begin{definition}
Let $X$ be a normal singular surface, $x\in $ Sing$(X),$ and $f:X^{\prime
}\rightarrow X$ a good resolution of $X.$ To the support $f^{-1}\left(
x\right) =\cup _{i=1}^{k}C_{i}$ of the exceptional divisor w.r.t. $f$
(resolving the singularity at $x)$ we can associate a \textit{weighted dual
graph} by assigning a weighted vertex to each $C_{i}$, with the weight being
the self intersection number $C_{i}^{2},$ and linking two vertices
corresponding to $C_{i}$ and $C_{j}$ by an edge of weight $(C_{i}\cdot
C_{j}).$ The \textit{fundamental cycle}
\begin{equation*}
Z_{\text{fund}}=\sum_{i=1}^{k}n_{i}C_{i},\ n_{i}>0,\ \ \forall i,\ \ 1\leq
i\leq k,
\end{equation*}
of $f$ w.r.t. $\left( X,x\right) $ is the unique, smallest positive cycle
for which $(Z_{\text{fund}}\cdot C_{i})\leq 0,$ for all $i,$ $1\leq i\leq k.$
\end{definition}

\begin{theorem}
\emph{(Artin \cite[Thm. 3]{ARTIN})} The following statements are equivalent%
\emph{:}\newline
\emph{(i) }$\left( X,x\right) $ is a rational surface singularity.\newline
\emph{(ii)} $p_{a}(Z_{\text{\emph{fund}}})=0.$ \emph{(}$p_{a}$ denotes here
the arithmetic genus\emph{).}
\end{theorem}

\begin{corollary}
\emph{(Brieskorn \cite[Lemma 1.3]{BRIESKORN2}) }For $\left( X,x\right) $ a
rational surface singularity, $\cup _{i=1}^{k}C_{i}$ has the following
properties\emph{:\smallskip\ }\newline
\emph{(i) }all $C_{i}$'s are smooth rational curves.\smallskip\ \newline
\emph{(ii)} $C_{i}\cap C_{j}\cap C_{l}=\varnothing $ for pairwise distinct $%
i,j,l.\smallskip $ \newline
\emph{(iii) }$(C_{i}\cdot C_{j})\in \left\{ 0,1\right\} ,$ for $i\neq j.$
\smallskip \newline
\emph{(iv)} The weighted dual graph contains no cycles.
\end{corollary}

\begin{corollary}
\emph{(Artin \cite[Cor. 6]{ARTIN})} \ If $\left( X,x\right) $ is a rational
surface singularity,
\begin{equation*}
\emph{mult}_{x}\left( X\right) =\emph{mult}\left( \mathcal{O}_{X,x}\right)
\end{equation*}
the \emph{multiplicity} of $X$ at $x$ and
\begin{equation*}
\emph{edim}\left( X,x\right) =\emph{dim}_{\mathbb{C}}\left( \frak{m}_{X,x}/%
\frak{m}_{X,x}^{2}\right)
\end{equation*}
its \emph{minimal embedding dimension}, then we have\emph{:}
\begin{equation*}
-Z_{\text{\emph{fund}}}^{2}=\text{\emph{mult}}_{x}\left( X\right) =\text{%
\emph{edim}}\left( X,x\right) -1.
\end{equation*}
\end{corollary}

\begin{theorem}
\label{RDP}The following conditions for a normal surface singularity $\left(
X,x\right) $ are equivalent\emph{:\smallskip }\newline
\emph{(i) }$\left( X,x\right) $ is a canonical singularity.\smallskip\
\newline
\emph{(ii)} $\left( X,x\right) $ is a rational Gorenstein
singularity.\smallskip\ \newline
\emph{(iii)} $\left( X,x\right) $ is a \emph{rational double point (RDP)} or
a \emph{Kleinian} or \emph{Du} \emph{Val singularity,} i.e, it is
analytically equivalent to the hypersurface singularity
\begin{equation*}
(\left\{ \left( z_{1},z_{2},z_{3}\right) \in \mathbb{C}^{3}\ \left| \
\varphi \left( z_{1},z_{2},z_{3}\right) =0\right. \right\} ,(0,0,0))
\end{equation*}
which is determined by one of the quasihomogeneous polynomials of type $A$-$%
D $-$E$ of the table\emph{:\medskip }
\begin{equation*}
\begin{tabular}{|l|l|}
\hline
$\mathbf{Type}$ & $\varphi \left( z_{1},z_{2},z_{3}\right) $ \\ \hline\hline
$\mathbf{A}_{n}\ \ \ (n\geq 1)$ & $
\begin{array}{c}
\, \\
z_{1}^{n+1}+z_{2}^{2}+z_{3}^{2} \\
\,
\end{array}
$ \\ \hline
$\mathbf{D}_{n}\ \ \ (n\geq 4)$ & $
\begin{array}{c}
\, \\
z_{1}^{n-1}+z_{1}z_{2}^{2}+z_{3}^{2} \\
\,
\end{array}
$ \\ \hline
$\mathbf{E}_{6}$ & $
\begin{array}{c}
\, \\
z_{1}^{4}+z_{2}^{3}+z_{3}^{2} \\
\,
\end{array}
$ \\ \hline
$\mathbf{E}_{7}$ & $
\begin{array}{c}
\, \\
z_{1}^{3}z_{2}+z_{2}^{3}+z_{3}^{2} \\
\,
\end{array}
$ \\ \hline
$\mathbf{E}_{8}$ & $
\begin{array}{c}
\, \\
z_{1}^{5}+z_{2}^{3}+z_{3}^{2} \\
\,
\end{array}
$ \\ \hline
\end{tabular}
\end{equation*}
\newline
\emph{(iv)} $(X,x)$ is analytically equivalent to a quotient singularity $%
\left( \mathbb{C}^{2}/G,[\mathbf{0}]\right) ,$ where $G$ denotes a finite
subgroup of \emph{SL}$\left( 2,\mathbb{C}\right) .$ More precisely, taking
into account the classification \emph{(}up to conjugacy\emph{)} of these
groups \emph{(see \cite{DE LA HARPE-SIEGFRIED}, \cite{DURFEE}, \cite[p. 35]
{LAMOTKE}, \cite[\S 4.4]{SPRINGER}),} we get the correspondence\emph{%
:\medskip \medskip }
\begin{equation*}
\begin{tabular}{|l|l|l|l|l|l|}
\hline
$
\begin{array}{c}
\, \\
\mathbf{Type\ (as\ above)} \\
\,
\end{array}
$ & $\mathbf{A}_{n}$ & $\mathbf{D}_{n}\ $ & $\mathbf{E}_{6}$ & $\mathbf{E}%
_{7}$ & $\mathbf{E}_{8}$ \\ \hline\hline
$
\begin{array}{c}
\, \\
\mathbf{Type\ of\ the\ acting\ group\ }G \\
\,
\end{array}
$ & $\mathbf{C}_{n}$ & $\mathbf{D}_{n}$ & $\mathbf{T}$ & $\mathbf{O}$ & $%
\mathbf{I}$ \\ \hline
\end{tabular}
\end{equation*}
\newline
\emph{(}By $\mathbf{C}_{n}$ we denote a cyclic group of order $n$, and by $%
\mathbf{D}_{n},\mathbf{T,O}$ and $\mathbf{I}$ the binary dihedral,
tetrahedral, octahedral and icosahedral subgroups of \emph{SL}$\left( 2,%
\mathbb{C}\right) $, having orders $4(n-2),$ $24,48$ and $120$, respectively%
\emph{).}\medskip\ \newline
\emph{(v)\ [Inductive criterion] \ }$(X,x)$ is an \emph{absolutely isolated
double point}, in the sense, that for any finite sequence
\begin{equation*}
\{\pi _{j-1}:X_{j}=\mathbf{Bl}_{\left\{ x_{j-1}\right\} }^{\text{\emph{red}}%
}\longrightarrow X_{j-1}\ \left| \ 1\leq j\leq l\right. \}
\end{equation*}
of blow-ups with closed \emph{(}reduced\emph{)} points as centers and $%
X_{0}=X,$ the only singular points of $X_{l}$ are isolated double points.
\emph{(}In particular, $(X,x)$ is a hypersurface double point whose normal
cone is either a \emph{(}not necessarily irreducible\emph{)} plane conic or
a double line\emph{).}
\end{theorem}

\begin{definition}
Let $(X,x)$ be a normal surface singularity. Assume that $(X,x)$ is elliptic%
\footnote{%
Laufer \cite{Laufer4} calls an elliptic singularity (in the above sense)
``minimally elliptic''.} and Gorenstein. We define the Laufer-Reid invariant
LRI$\left( X,x\right) $ of $x$ in $X$ to be the self-intersection number of
Artin's fundamental cycle with opposite sign:
\begin{equation*}
\text{LRI}\left( X,x\right) =-Z_{\text{fund}}^{2}
\end{equation*}
\end{definition}

\begin{theorem}
\label{LAUFER-REID}\emph{(Laufer \cite{Laufer4}, Reid \cite{REID0})} Let $%
(X,x)$ be a normal surface singularity. Assume that $(X,x)$ is elliptic and
Gorenstein. Then \emph{LRI}$\left( X,x\right) \geq 1$ and has the following
properties\emph{:\medskip\ }\newline
\emph{(i) }If \emph{LRI}$\left( X,x\right) =1,$ then
\begin{equation*}
(X,x)\cong (\left\{ \left( z_{1},z_{2},z_{3}\right) \in \mathbb{C}^{3}\
\left| \ z_{1}^{2}+z_{2}^{3}+\varphi \left( z_{2},z_{3}\right) =0\right.
\right\} ,(0,0,0)),
\end{equation*}
with $\varphi \left( z_{2},z_{3}\right) $ a sum of monomials of the form $%
z_{2}z_{3}^{\kappa },$ $\kappa \in \mathbb{Z}_{\geq 4},$ and $z_{3}^{\kappa
},$ $\kappa \in \mathbb{Z}_{\geq 6}.$ Performing the monoidal transformation
\begin{equation*}
\mathbf{Bl}_{\left\{ x\right\} }^{\left( z_{1}^{3},z_{2}^{2},z_{3}\right)
}\longrightarrow X
\end{equation*}
we get a normal surface having at most one Du Val point.\medskip\ \newline
\emph{(ii)} If \emph{LRI}$\left( X,x\right) =2,$ then
\begin{equation*}
(X,x)\cong (\left\{ \left( z_{1},z_{2},z_{3}\right) \in \mathbb{C}^{3}\
\left| \ z_{1}^{2}+\varphi \left( z_{2},z_{3}\right) =0\right. \right\}
,(0,0,0)),
\end{equation*}
with $\varphi \left( z_{2},z_{3}\right) $ a sum of monomials of the form $%
z_{2}^{\kappa }z_{3}^{\lambda },\kappa +\lambda \in \mathbb{Z}_{\geq 4}.$ In
this case, the normalized blow-up $X$ at $x$%
\begin{equation*}
\text{\emph{Norm}}\left[ \mathbf{Bl}_{\left\{ x\right\} }^{\text{\emph{red}}}%
\right] =\mathbf{Bl}_{\left\{ x\right\} }^{\left(
z_{1}^{2},z_{2},z_{3}\right) }\longrightarrow X
\end{equation*}
has only Du Val singular points.\medskip\ \newline
$\emph{(iii)}$ If \emph{LRI}$\left( X,x\right) \geq 2,$ then
\begin{equation*}
\emph{LRI}\left( X,x\right) =\emph{mult}_{x}\left( X\right) .\smallskip
\newline
\end{equation*}
\emph{(iv)} If \emph{LRI}$\left( X,x\right) \geq 3,$ then
\begin{equation*}
\emph{LRI}\left( X,x\right) =\emph{edim}(X,x)
\end{equation*}
and the overlying space of the ordinary blow-up
\begin{equation*}
\mathbf{Bl}_{\left\{ x\right\} }^{\text{\emph{red}}}=\mathbf{Bl}_{\left\{
x\right\} }^{\frak{m}_{X,x}}\longrightarrow X
\end{equation*}
of $x$ is a normal surface with at most Du Val singularities.
\end{theorem}

\begin{theorem}
Let $(X,x)$ be a normal surface singularity. Then\emph{\footnote{%
Explanation of terminology: Let $f:X^{\prime }\rightarrow X$ be the good
minimal resolution of $X.$ Then $x$ is called \textit{simple-elliptic} (resp.%
\textit{\ a cusp}) if the support of the exceptional divisor w.r.t. $f$
lying over $x$ consists of a smooth elliptic curve (resp. a cycle of $%
\mathbb{P}_{\mathbb{C}}^{1}$'s).}\medskip \medskip \medskip }{\small
\begin{equation*}
\left\{
\begin{array}{cccc}
\text{\emph{(i)}} & x\text{ is \emph{terminal}} & \Longleftrightarrow & x\in
\text{\emph{Reg}}\left( X\right) \\
\  &  & \, &  \\
\text{\emph{(ii)}} & x\text{ is \emph{canonical}} & \Longleftrightarrow &
\left(
\begin{array}{c}
\left( X,x\right) \cong \left( \mathbb{C}^{2}/G,[\mathbf{0}]\right) \ \text{%
\emph{with}} \\
G\ \text{\emph{a finite subgroup of} \emph{SL}}(2,\mathbb{C})
\end{array}
\right) \\
\  &  & \, &  \\
\text{\emph{(iii)}} & x\text{ is \emph{log-terminal}} & \Longleftrightarrow
& \left(
\begin{array}{c}
\left( X,x\right) \cong \left( \mathbb{C}^{2}/G,[\mathbf{0}]\right) \ \text{%
\emph{with}} \\
G\ \text{\emph{a finite subgroup of} \emph{GL}}(2,\mathbb{C})
\end{array}
\right) \\
\  &  & \, &  \\
\text{\emph{(iv)}} & x\text{ is \emph{log-canonical}} & \Longleftrightarrow
& \left(
\begin{array}{c}
x\ \text{\emph{is simple-elliptic, a cusp }} \\
\text{\emph{or a regular point }} \\
\text{\emph{or a quotient thereof}}
\end{array}
\right)
\end{array}
\right.
\end{equation*}
}
\end{theorem}

\noindent {}Log-terminal surface singularities are rational (by Theorem \ref
{LOG-RAT}). This is, of course, not the case for log-canonical surface
singularities which are not log-terminal. For the fine classification of
log-terminal and log-canonical surface singularities, the reader is referred
to the papers of Brieskorn\footnote{%
Brieskorn classified (up to conjugation) \textit{all} finite subgroups of GL$%
(2,\mathbb{C)}$ in \cite[2.10 and 2.11]{BRIESKORN2}.} \cite{BRIESKORN2},
Iliev \cite{ILIEV}, Kawamata \cite{KAWAMATA}, Alexeev \cite{ALEXEEV} and
Blache \cite{BLACHE}. \bigskip

\noindent {}$\bullet $ Next, let us recall some basic facts from the theory
of 3-dimensional terminal and canonical singularities.

\begin{definition}
A normal threefold singularity $(X,x)$ is called \textit{compound Du Val
singularity} (abbreviated: \textit{cDV singularity}) if for some general
hyperplane section $\mathbb{H}$ through $x,$ $(\mathbb{H},x)$ is a Du Val
singularity, or equivalently, if{\small
\begin{equation*}
(X,x)\cong (\{\left( z_{1},z_{2},z_{3},z_{4}\right) \in \mathbb{C}%
^{4}\,\left| \,\varphi \left( z_{1},z_{2},z_{3}\right) +z_{4}\cdot g\left(
z_{1},z_{2},z_{3},z_{4}\right) =0\right. \},\left( 0,0,0,0\right) ),
\end{equation*}
}where $\varphi \left( z_{1},z_{2},z_{3}\right) $ is one of the
quasihomogeneous polynomials listed in the Thm. \ref{RDP} (iii) and $g\left(
z_{1},z_{2},z_{3},z_{4}\right) $ an arbitrary polynomial in $\mathbb{C}%
[z_{1},z_{2},z_{3},z_{4}].$ According to the type of $\varphi \left(
z_{1},z_{2},z_{3}\right) ,(X,x)$ is called $c\mathbf{A}_{n},c\mathbf{D}_{n},c%
\mathbf{E}_{6},c\mathbf{E}_{7}$ and $c\mathbf{E}_{8}$-point, respectively.
Compound Du Val singularities are not necessarily isolated.
\end{definition}

\begin{theorem}
\emph{(Reid \cite[0.6 (I), 1.1, 1.11]{REID2}) }Let $(X,x)$ be a normal
threefold singularity. Then\emph{\medskip }{\small \medskip\
\begin{equation*}
\left\{
\begin{array}{cccc}
\text{\emph{(i)}} & x\ \text{\emph{is terminal}} & \Longrightarrow & x\
\text{\emph{is isolated}} \\
&  & \, &  \\
\text{\emph{(ii)}} & x\ \text{\emph{is terminal of index} }1 &
\Longleftrightarrow & x\ \text{\emph{is an isolated cDV point}} \\
&  & \, &  \\
\text{\emph{(iii)}} & x\ \text{\emph{is terminal of index }}\geq 1\text{%
\emph{\ }} & \Longleftrightarrow & \left(
\begin{array}{c}
x\ \text{\emph{is a quotient of an isolated}} \\
\text{\emph{cDV point by a finite cyclic group}}
\end{array}
\right) \\
&  & \, &  \\
\text{\emph{(iv)}} & x\ \text{\emph{is a cDV point}} & \Longrightarrow & x\
\text{\emph{is canonical}}
\end{array}
\right.
\end{equation*}
}
\end{theorem}

\noindent{}For the extended lists of the fine classification of
3-dimensional terminal singularities of arbitrary index, see Mori \cite{MORI}%
, Reid \cite{REID4} and Koll\'{a}r \& Shepherd-Barron \cite{KO-SB}. The
normal forms of the defining equations of cDV points have been studied by
Markushevich in \cite{MARKUSHEVICH}. On the other hand, 3-dimensional
terminal cyclic quotient singularities, which play a crucial role in the
above cited investigations, are quite simple.

\begin{theorem}
\emph{(Danilov \cite{DANILOV}, Morrison-Stevens \cite{Mo-Ste})} Let $(X,x)$
be a terminal threefold singularity. Then\emph{\medskip }{\small
\begin{equation*}
\begin{array}{c}
\left(
\begin{array}{c}
(X,x)\cong \left( \mathbb{C}^{3}/G,[\mathbf{0}]\right) \ \text{\emph{with} }G%
\emph{\ }\text{\emph{a linearly acting finite}} \\
\text{\emph{cyclic subgroup of GL}}(3,\mathbb{C})\ \text{\emph{without
pseudoreflections}}
\end{array}
\right) \\
\, \\
\Updownarrow \\
\, \\
\left(
\begin{array}{c}
\text{\emph{the action of} }G\ \text{\emph{is given (up to permutations }}
\\
\text{\emph{of }}(z_{1},z_{2},z_{3})\ \text{\emph{and group symmetries) by}}
\\
(z_{1},z_{2},z_{3})\longmapsto (\zeta _{\mu }^{\lambda }\ z_{1},\zeta _{\mu
}^{-\lambda }\ z_{2},\zeta _{\mu }\ z_{3}),\ \text{\emph{where} }\mu
:=\left| G\right| , \\
\gcd \left( \lambda ,\mu \right) =1,\ \text{\emph{and} }\zeta _{\mu }\ \text{%
\emph{denotes a} }\mu \text{\emph{-th root of unity}}
\end{array}
\right) \\
\,
\end{array}
\end{equation*}
}
\end{theorem}

\noindent {}$\bullet $ \textbf{Reduction of 3-dimensional canonical
singularities}. The singularities of a quasiprojective threefold $X$ can be
reduced by a ``canonical modification'' $X^{\text{can}}\overset{g}{%
\rightarrow }X,$ so that $K_{X^{\text{can}}}$ is $g$-ample. $X^{\text{can}}$
can be also modified by a ``terminal modification'' $X^{\text{ter}}\overset{f%
}{\rightarrow }X^{\text{can}},$ so that $X^{\text{ter}}$ has at most
terminal singularities, where $f$ is projective and crepant. Finally, $X^{%
\text{ter}}$ can be modified by another modification $X^{\mathbb{Q}\text{%
-f-ter}}\overset{h}{\rightarrow }X^{\text{ter}}$, so that $X^{\mathbb{Q}%
\text{-f-ter}}$ has at most $\mathbb{Q}$-factorial terminal singularities%
\footnote{%
The morphism $h$ can be constructed by taking successively birational
morphisms of the form
\begin{equation*}
\mathbf{Proj}(\bigoplus\limits_{\nu \geq 0}\mathcal{O}_{X^{\text{ter}%
}}\left( \nu D\right) )\longrightarrow X^{\text{ter}},
\end{equation*}
where $D$'s are Weil divisors which are not $\mathbb{Q}$-Cartier divisors
(cf. \cite[p. 201]{Ko-Mo}).}, and $h$ is projective and an isomorphism in
codimension $1$. (See \cite{REID1}, \cite{REID2}, \cite{REID4}, \cite
{KAWAMATA}, \cite{MORI2}, and \cite[section 6.3]{Ko-Mo}). The main steps of
the intrinsic construction of $f$, due to Miles Reid, will be explained in
broad outline and will be applied in the framework of toric geometry in
section \ref{SECTION4}.\medskip

\noindent {}\underline{\textbf{Step 1.}}\textbf{\ \ Reduction to index 1
canonical singularities by index cover.} If $x\in X:=X^{\text{can}}$ is a
canonical singularity of index $\ell >1,$ then one considers the finite
Galois cover
\begin{equation*}
\phi :Y=\text{Spec}\left( \bigoplus\limits_{i=0}^{\ell -1}\mathcal{O}%
_{X}\left( \ell K_{X}\right) \right) \longrightarrow X.
\end{equation*}
The preimage $\phi ^{-1}\left( x\right) $ constists of just one point, say $%
y,$ and if $y\in Y$ is terminal, then the same is also valid for $x\in X.$
Moreover, if $\psi :Y^{\prime }\rightarrow Y$ is a crepant resolution of $Y$
(as those ones which will be constructed in the next steps), then we get a
commutative diagram
\begin{equation*}
\begin{array}{ccc}
Y^{\prime } & \overset{\psi }{\longrightarrow } & Y \\
\downarrow \phi ^{\prime } &  & \phi \downarrow \\
Y^{\prime }/\mathbb{Z}_{\ell } & \overset{\psi ^{\prime }}{\longrightarrow }
& X
\end{array}
\end{equation*}
extending the action of $\mathbb{Z}_{\ell }$ on $\left( X,x\right) $ to an
action on $Y^{\prime }$ where $\psi ^{\prime }$ is crepant with at least one
exceptional prime divisor and $\phi ^{\prime }$ is etale in codimension $1$.
\medskip

\noindent {}\underline{\textbf{Step 2.}} \textbf{Weighted blow-ups of
non-cDv singularities.} From now on we may assume that $X$ contains at most
canonical singularities of index $1$ (i.e., rational Gorenstein
singularities). If $\ X$ contains non-cDV points $x\in X$, then for a
general hyperplane section $\mathbb{H}$ through $x,$ $(\mathbb{H},x)$ is an
elliptic surface singularity. Using Theorem \ref{LAUFER-REID} one obtains
the following:

\begin{proposition}
\label{PROP-REID}\emph{(i) }If \emph{LRI}$\left( \mathbb{H},x\right) =1,$
then
\begin{equation*}
(X,x)\cong (\left\{ \left( z_{1},z_{2},z_{3},z_{4}\right) \in \mathbb{C}%
^{4}\ \left| \ z_{1}^{2}+z_{2}^{3}+\varphi \left( z_{2},z_{3},z_{4}\right)
=0\right. \right\} ,(0,0,0,0)),
\end{equation*}
with $\varphi \left( z_{2},z_{3},z_{4}\right)
=z_{2}F_{1}(z_{3},z_{4})+F_{2}(z_{3},z_{4}),$ where $F_{1}$ \emph{(}resp. $%
F_{2}$\emph{)} is a sum of monomials $z_{3}^{\kappa }z_{4}^{\lambda }$ of
degree $\kappa +\lambda \geq 4$ \emph{(}resp. $\geq 6$\emph{)}.\smallskip
\newline
\emph{(ii)} If \emph{LRI}$\left( \mathbb{H},x\right) =2,$ then
\begin{equation*}
(X,x)\cong (\left\{ \left( z_{1},z_{2},z_{3},z_{4}\right) \in \mathbb{C}%
^{4}\ \left| \ z_{1}^{2}+\varphi \left( z_{2},z_{3},z_{4}\right) =0\right.
\right\} ,(0,0,0,0)),
\end{equation*}
with $\varphi \left( z_{2},z_{3},z_{4}\right) $ a sum of monomials of degree
$\geq 4.\smallskip $\newline
\emph{(iii)} If \emph{LRI}$\left( \mathbb{H},x\right) \geq 3,$ then \emph{LRI%
}$\left( \mathbb{H},x\right) =\emph{edim}\left( \mathbb{H},x\right) =$ \emph{%
edim}$(X,x)-1.$
\end{proposition}

\noindent {}Blowing up $x\in X$ with respect to the weights $\left(
2,1,1,1\right) ,\left( 3,2,1,1\right) $ and $(1,1,1,1)$ for LRI$\left(
\mathbb{H},x\right) =1,2$ and $\geq 3,$ respectively, we get a projective
crepant partial desingularization of $X.$ Repeating this procedure for all
the non-cDV points of $X,$ we reduce our singularities to cDV
singularities.\medskip

\noindent {}\underline{\textbf{Step 3.}} \textbf{Simultaneous blow-up of
one-dimensional singular loci.} From now on we may assume that $X$ contains
at most cDV singularities. If Sing$\left( X\right) $ contains
one-dimensional components, then we blow their union up (by endowing it with
the reduced subscheme structure). This blow-up is realized by a projective,
crepant birational morphism. Repeating this procedure finitely many times we
reduce our singularities to isolated cDV singularities, i.e., to terminal
singularities of index $1$.

\begin{remark}
After step 3, one may use the above projective birational morphism $h$ to
get only $\mathbb{Q}$-factorial terminal singularities. Sometimes, it is
also useful to desingularize overall our threefold by resolving the
remaining non-$\mathbb{Q}$-factorial terminal singularities.
\end{remark}

\section{Toric singularities\label{SECTION3}}

\noindent Toric singularities occupy a distinguished position within the
class of rational singularities, as they can be described by binomial-type
equations. In this section we shall introduce the brief toric glossary
\textsf{(a)}-\textsf{(k) }and the notation which will be used in the sequel,
and we shall summarize their main properties. For further details on toric
geometry the reader is referred to the textbooks of Oda \cite{Oda}, Fulton
\cite{Fulton} and Ewald \cite{Ewald}, and to the lecture notes \cite{KKMS}.
\medskip

\noindent \textsf{(a)} The \textit{linear hull, }the\textit{\ affine hull},
the \textit{positive hull} and \textit{the convex hull} of a set $B$ of
vectors of $\mathbb{R}^{r}$, $r\geq 1,$ will be denoted by lin$\left(
B\right) $, aff$\left( B\right) $, pos$\left( B\right) $ (or $\mathbb{R}%
_{\geq 0}\,B$) and conv$\left( B\right) $, respectively. The \textit{%
dimension} dim$\left( B\right) $ of a $B\subset \mathbb{R}^{r}$ is defined
to be the dimension of aff$\left( B\right) $. \newline
\newline
\textsf{(b) }Let $N$ be a free $\mathbb{Z}$-module of rank $r\geq 1$. $N$
can be regarded as a \textit{lattice }in $N_{\mathbb{R}}:=N\otimes _{\mathbb{%
Z}}\mathbb{R}\cong \mathbb{R}^{r}$. The \textit{lattice determinant} det$%
\left( N\right) $ of $N$ is the $r$-volume of the parallelepiped spanned by
any $\mathbb{Z}$-basis of it. An $n\in N$ is called \textit{primitive} if
conv$\left( \left\{ \mathbf{0},n\right\} \right) \cap N$ contains no other
points except $\mathbf{0}$ and $n$.\smallskip

Let $N$ be as above, $M:=$ Hom$_{\mathbb{Z}}\left( N,\mathbb{Z}\right) $ its
dual lattice, $N_{\mathbb{R}},M_{\mathbb{R}}$ their real scalar extensions,
and $\left\langle .,.\right\rangle :M_{\mathbb{R}}\times N_{\mathbb{R}%
}\rightarrow \mathbb{R}$ the natural $\mathbb{R}$-bilinear pairing. A subset
$\sigma $ of $N_{\mathbb{R}}$ is called \textit{convex polyhedral cone} (%
\textit{c.p.c.}, for short) if there exist $n_{1},\ldots ,n_{k}\in N_{%
\mathbb{R}}$, such that
\begin{equation*}
\sigma =\text{pos}\left( \left\{ n_{1},\ldots ,n_{k}\right\} \right) \,.
\end{equation*}
Its \textit{relative interior }int$\left( \sigma \right) $ is the usual
topological interior of it, considered as subset of lin$\left( \sigma
\right) =\sigma +\left( -\sigma \right) $. The \textit{dual cone} $\sigma
^{\vee }$ of a c.p.c. $\sigma $ is a c.p. cone defined by
\begin{equation*}
\sigma ^{\vee }:=\left\{ \mathbf{y}\in M_{\mathbb{R}}\ \left| \ \left\langle
\mathbf{y},\mathbf{x}\right\rangle \geq 0,\ \forall \mathbf{x},\ \mathbf{x}%
\in \sigma \right. \right\} \;.\;
\end{equation*}
Note that $\left( \sigma ^{\vee }\right) ^{\vee }=\sigma $ and
\begin{equation*}
\text{dim}\left( \sigma \cap \left( -\sigma \right) \right) +\text{dim}%
\left( \sigma ^{\vee }\right) =\text{dim}\left( \sigma ^{\vee }\cap \left(
-\sigma ^{\vee }\right) \right) +\text{dim}\left( \sigma \right) =r.
\end{equation*}
A subset $\tau $ of a c.p.c. $\sigma $ is called a \textit{face} of $\sigma $
(notation: $\tau \prec \sigma $), if
\begin{equation*}
\tau =\left\{ \mathbf{x}\in \sigma \ \left| \ \left\langle m_{0},\mathbf{x}%
\right\rangle =0\right. \right\} ,
\end{equation*}
for some $m_{0}\in \sigma ^{\vee }$. A c.p.c. $\sigma =$ pos$\left( \left\{
n_{1},\ldots ,n_{k}\right\} \right) $ is called \textit{simplicial} (resp.
\textit{rational}) if $n_{1},\ldots ,n_{k}$ are $\mathbb{R}$-linearly
independent (resp. if $n_{1},\ldots ,n_{k}\in N_{\mathbb{Q}}$, where $N_{%
\mathbb{Q}}:=N\otimes _{\mathbb{Z}}\mathbb{Q}$). A \textit{strongly convex
polyhedral cone }(\textit{s.c.p.c.}, for short) is a c.p.c. $\sigma $ for
which $\sigma \cap \left( -\sigma \right) =\left\{ \mathbf{0}\right\} $,
i.e., for which dim$\left( \sigma ^{\vee }\right) =r$. The s.c.p. cones are
alternatively called \textit{pointed cones} (having $\mathbf{0}$ as their
apex).\newline
\newline
\textsf{(c) }If $\sigma \subset N_{\mathbb{R}}$ is a c.p. cone, then the
subsemigroup $\sigma \cap N$ of $N$ is a monoid. The following proposition
is due to Gordan, Hilbert and van der Corput and describes its fundamental
properties.

\begin{proposition}[Minimal generating system]
\label{MINGS}If $\sigma \subset N_{\mathbb{R}}$ is a c.p. rational cone,
then $\sigma \cap N$ is finitely generated as additive semigroup. Moreover,
if $\sigma $ is strongly convex, then among all the systems of generators of
$\sigma \cap N$, there is a system $\mathbf{Hilb}_{N}\left( \sigma \right) $
of \emph{minimal cardinality}, which is uniquely determined \emph{(}up to
the ordering of its elements\emph{)} by the following characterization\emph{%
:\smallskip }
\begin{equation}
\mathbf{Hilb}_{N}\left( \sigma \right) =\left\{ n\in \sigma \cap \left(
N\smallsetminus \left\{ \mathbf{0}\right\} \right) \ \left| \
\begin{array}{l}
n\ \text{\emph{cannot be expressed }} \\
\text{\emph{as the sum of two other vectors}} \\
\text{\emph{belonging to\ } }\sigma \cap \left( N\smallsetminus \left\{
\mathbf{0}\right\} \right)
\end{array}
\right. \right\}  \label{Hilbbasis}
\end{equation}
$\mathbf{Hilb}_{N}\left( \sigma \right) $ \emph{is called }\textit{the
Hilbert basis of }$\sigma $ w.r.t. $N.\medskip $
\end{proposition}

\noindent \textsf{(d)} For a lattice $N$ of rank $r$ having $M$ as its dual,
we define an $r$-dimensional \textit{algebraic torus }$T_{N}\cong \left(
\mathbb{C}^{\ast }\right) ^{r}$ by setting $T_{N}:=$ Hom$_{\mathbb{Z}}\left(
M,\mathbb{C}^{\ast }\right) =N\otimes _{\mathbb{Z}}\mathbb{C}^{\ast }$.
Every $m\in M$ assigns a character $\mathbf{e}\left( m\right)
:T_{N}\rightarrow \mathbb{C}^{\ast }$. Moreover, each $n\in N$ determines a
$1$-parameter subgroup
\begin{equation*}
\vartheta _{n}:\mathbb{C}^{\ast }\rightarrow T_{N}\ \ \ \text{with\ \ \ }%
\vartheta _{n}\left( \lambda \right) \left( m\right) :=\lambda
^{\left\langle m,n\right\rangle }\text{, \ \ for\ \ \ }\lambda \in \mathbb{C}%
^{\ast },\ m\in M\ .\
\end{equation*}
We can therefore identify $M$ with the character group of $T_{N}$ and $N$
with the group of $1$-parameter subgroups of $T_{N}$. On the other hand, for
a rational s.c.p.c. $\sigma $ with $M\cap \sigma ^{\vee }=\mathbb{Z}_{\geq
0}\ m_{1}+\cdots +\mathbb{Z}_{\geq 0}\ m_{\nu }$, we associate to the
finitely generated monoidal subalgebra
\begin{equation*}
\mathbb{C}\left[ M\cap \sigma ^{\vee }\right] =\bigoplus_{m\in M\cap \sigma
^{\vee }}\mathbf{e}\left( m\right)
\end{equation*}
of the $\mathbb{C}$-algebra $\mathbb{C}\left[ M\right] =\oplus _{m\in M}%
\mathbf{e}\left( m\right) $ an affine complex variety\footnote{%
As point-set $U_{\sigma }$ is actually the ``maximal spectrum'' Max-Spec$%
\left( \mathbb{C}\left[ M\cap \sigma ^{\vee }\right] \right) .$}
\begin{equation*}
U_{\sigma }:=\text{Spec}\left( \mathbb{C}\left[ M\cap \sigma ^{\vee }\right]
\right) ,
\end{equation*}
which can be identified with the set of semigroup homomorphisms :
\begin{equation*}
U_{\sigma }=\left\{ u:M\cap \sigma ^{\vee }\ \rightarrow \mathbb{C\ }\left|
\begin{array}{c}
\ u\left( \mathbf{0}\right) =1,\ u\left( m+m^{\prime }\right) =u\left(
m\right) \cdot u\left( m^{\prime }\right) ,\smallskip \  \\
\text{for all \ \ }m,m^{\prime }\in M\cap \sigma ^{\vee }
\end{array}
\right. \right\} \ ,
\end{equation*}
where $\mathbf{e}\left( m\right) \left( u\right) :=u\left( m\right) ,\
\forall m,\ m\in M\cap \sigma ^{\vee }\ $ and\ $\forall u,\ u\in U_{\sigma }$%
.

\begin{proposition}[Embedding by binomials]
\label{EMB}In the analytic category, $U_{\sigma }$, identified with its
image under the injective map
\begin{equation*}
\left( \mathbf{e}\left( m_{1}\right) ,\ldots ,\mathbf{e}\left( m_{\nu
}\right) \right) :U_{\sigma }\hookrightarrow \mathbb{C}^{\nu },
\end{equation*}
can be regarded as an analytic set determined by a system of equations of
the form\emph{:} \emph{(monomial) = (monomial).} This analytic structure
induced on $U_{\sigma }$ is independent of the semigroup generators $\left\{
m_{1},\ldots ,m_{\nu }\right\} $ and each map $\mathbf{e}\left( m\right) $
on $U_{\sigma }$ is holomorphic w.r.t. it. In particular, for $\tau \prec
\sigma $, $U_{\tau }$ is an open subset of $U_{\sigma }$. Moreover, if $%
\#\left( \mathbf{Hilb}_{M}\left( \sigma ^{\vee }\right) \right) =k\ \left(
\leq \nu \right) $, then, by \emph{(\ref{Hilbbasis}), }$k$\emph{\ }is
nothing but the \emph{embedding dimension} of $U_{\sigma }$, i.e., the \emph{%
minimal} number of generators of the maximal ideal of the local $\mathbb{C}$%
-algebra $\mathcal{O}_{U_{\sigma },\ \mathbf{0}}^{\text{\emph{hol}}}$.
\end{proposition}

\noindent \textit{Proof. }See Oda \cite[Prop. 1.2 and 1.3., pp. 4-7]{Oda}%
.\hfill$\square \bigskip $

\noindent \textsf{(e)} A \textit{fan }w.r.t. a free $\mathbb{Z}$-module%
\textit{\ }$N$ is a finite collection $\Delta $ of rational s.c.p. cones in $%
N_{\mathbb{R}}$, such that :\smallskip \newline
(i) any face $\tau $ of $\sigma \in \Delta $ belongs to $\Delta $,
and\smallskip \newline
(ii) for $\sigma _{1},\sigma _{2}\in \Delta $, the intersection $\sigma
_{1}\cap \sigma _{2}$ is a face of both $\sigma _{1}$ and $\sigma
_{2}.\smallskip $\newline
By $\left| \Delta \right| :=\cup \left\{ \sigma \ \left| \ \sigma \in \Delta
\right. \right\} $ one denotes the \textit{support} and by $\Delta \left(
i\right) $ the set of all $i$-dimensional cones of a fan $\Delta $ for $%
0\leq i\leq r$. If $\varrho \in \Delta \left( 1\right) $ is a ray, then
there exists a unique primitive vector $n\left( \varrho \right) \in N\cap
\varrho $ with $\varrho =\mathbb{R}_{\geq 0}\ n\left( \varrho \right) $ and
each cone $\sigma \in \Delta $ can be therefore written as
\begin{equation*}
\sigma =\sum_{\varrho \in \Delta \left( 1\right) ,\ \varrho \prec \sigma }\
\mathbb{R}_{\geq 0}\ n\left( \varrho \right) \ \ .
\end{equation*}
The set
\begin{equation*}
\text{Gen}\left( \sigma \right) :=\left\{ n\left( \varrho \right) \ \left| \
\varrho \in \Delta \left( 1\right) ,\varrho \prec \sigma \right. \right\}
\end{equation*}
is called the\textit{\ set of minimal generators }(within the pure first
skeleton) of $\sigma $. For $\Delta $ itself one defines analogously Gen$%
\left( \Delta \right) :=\bigcup_{\sigma \in \Delta }$ Gen$\left( \sigma
\right) .\bigskip $\newline
\textsf{(f) }The \textit{toric variety X}$\left( N,\Delta \right) $
associated to a fan\textit{\ }$\Delta $ w.r.t. the lattice\textit{\ }$N$ is
by definition the identification space
\begin{equation}
X\left( N,\Delta \right) :=((\coprod\limits_{\sigma \in \Delta }\ U_{\sigma
})\ /\ \sim )  \label{TORVAR}
\end{equation}
with $U_{\sigma _{1}}\ni u_{1}\sim u_{2}\in U_{\sigma _{2}}$ if and only if
there is a $\tau \in \Delta ,$ such that $\tau \prec \sigma _{1}\cap \sigma
_{2}$ and $u_{1}=u_{2}$ within $U_{\tau }$. $X\left( N,\Delta \right) $ is
called \textit{simplicial} if all the cones of $\Delta $ are simplicial. $%
X\left( N,\Delta \right) $ is compact iff $\left| \Delta \right| =N_{\mathbb{%
R}}$ (\cite{Oda}, thm. 1.11, p. 16). Moreover, $X\left( N,\Delta \right) $
admits a canonical $T_{N}$-action which extends the group multiplication of $%
T_{N}=U_{\left\{ \mathbf{0}\right\} }$:
\begin{equation}
T_{N}\times X\left( N,\Delta \right) \ni \left( t,u\right) \longmapsto
t\cdot u\in X\left( N,\Delta \right)  \label{torus action}
\end{equation}
where, for $u\in U_{\sigma }\subset X\left( N,\Delta \right) $,
\begin{equation*}
\left( t\cdot u\right) \left( m\right) :=t\left( m\right) \cdot u\left(
m\right) ,\ \forall m,\ m\in M\cap \sigma ^{\vee }.
\end{equation*}
The orbits w.r.t. the action (\ref{torus action}) are parametrized by the
set of all the cones belonging to $\Delta $. For a $\tau \in \Delta $, we
denote by orb$\left( \tau \right) $ (resp. by $V\left( \tau \right) $) the
orbit (resp. the closure of the orbit) which is associated to $\tau $%
.\bigskip

\noindent \textsf{(g) }The group of $T_{N}$-invariant Weil divisors of a
toric variety $X\left( N,\Delta \right) $ has the set $\{V(\varrho )\,\left|
\,\varrho \in \Delta (1)\right. \}$ as $\mathbb{Z}$-basis. In fact, such a
divisor $D$ is of the form $D=D_{\psi }$, where
\begin{equation*}
D_{\psi }:=-\sum_{\varrho \in \Delta \left( 1\right) }\psi (n(\varrho
))V(\varrho )
\end{equation*}
and $\psi :\left| \Delta \right| \rightarrow \mathbb{R}$ a \textit{PL-}$%
\Delta $\textit{-support function}, i.e., an $\mathbb{R}$-valued, positively
homogeneous function on $\left| \Delta \right| $ with $\psi (N\cap \left|
\Delta \right| )\subset \mathbb{Z}$ which is piecewise linear and upper
convex on each $\sigma \in \Delta $. (\textit{Upper convex} on a $\sigma \in
\Delta $ means that $\psi \left| _{\sigma }(\mathbf{x}+\mathbf{x}^{\prime
})\right. \geq \psi \left| _{\sigma }(\mathbf{x})+\right. \psi \left|
_{\sigma }(\mathbf{x}^{\prime })\right. $, for all $\mathbf{x},\mathbf{x}%
^{\prime }\in \sigma $). For example, the canonical divisor $K_{X\left(
N,\Delta \right) }$ of $X\left( N,\Delta \right) $ equals $D_{\psi }$ for $%
\psi $ a PL-$\Delta $-support function\textit{\ }with $\psi (n(\varrho ))=1$%
, for all rays $\varrho \in \Delta \left( 1\right) $. A divisor $D=D_{\psi }$
is Cartier iff $\psi $ is a \textit{linear} $\Delta $-support function
(i.e., $\psi \left| _{\sigma }\right. $ is overall linear on each $\sigma
\in \Delta $). Obviously, $D_{\psi }$ is $\mathbb{Q}$-Cartier iff $k\cdot
\psi $ is a linear $\Delta $-support function for some $k\in \mathbb{N}$.

\begin{theorem}[Ampleness criterion]
\label{AMPLE}A $T_{N}$-invariant $\mathbb{Q}$-Cartier divisor $D=D_{\psi }$
of a toric variety $X\left( N,\Delta \right) $ of dimension $r$ is ample if
and only if there exists a $\kappa \in \mathbb{N}$, such that $\kappa \cdot
\psi $ is a \emph{strictly upper convex} linear $\Delta $-support function,
i.e., iff \ for every $\sigma \in \Delta (r)$ there is a unique $m_{\sigma
}\in M=$ \emph{Hom}$_{\mathbb{Z}}(N,\mathbb{Z})$, such that
\begin{equation*}
\kappa \cdot \psi (\mathbf{x})\leq \left\langle m_{\sigma },\mathbf{x}%
\right\rangle ,\text{ \ for all }\mathbf{x}\in \left| \Delta \right| ,
\end{equation*}
with equality being valid iff $\mathbf{x}\in \sigma $.
\end{theorem}

\noindent \textit{Proof. }It follows from \cite[Thm. 13, p. 48]{KKMS}%
.\smallskip \smallskip \hfill $\square \medskip $

\noindent {}\textsf{(h)}\textit{\ }The behaviour of toric varieties with
regard to the algebraic properties (\ref{HIERARCHY}) is as follows.

\begin{theorem}[Normality and CM-property]
All toric varieties are normal and Cohen-Macaulay.
\end{theorem}

\noindent \textit{Proof. }For a proof of the normality property see
\cite[Thm. 1.4, p. 7]{Oda}. The CM-property for toric varieties was first
shown by Hochster in \cite{Hochster}. See also Kempf \cite[Thm. 14, p. 52]
{KKMS}, and Oda \cite[3.9, p. 125]{Oda}.\hfill$\square \bigskip $\newline
In fact, by the definition (\ref{TORVAR}) of $X\left( N,\Delta \right) $,
all the algebraic properties of this kind are \textit{local }with respect to
its affine covering, i.e., it is enough to be checked for the affine toric
varieties $U_{\sigma }$ for all (maximal) cones $\sigma $ of the fan $\Delta
$.

\begin{definition}[Multiplicities and basic cones]
Let\emph{\ }$N$\emph{\ }be a free\emph{\ }$\mathbb{Z}$\emph{-}module of rank%
\emph{\ }$r$\emph{\ }and $\sigma \subset N_{\mathbb{R}}$ a simplicial,
rational s.c.p.c.\emph{\ }of dimension $d\leq r$. $\sigma $ can be obviously
written as $\sigma =\varrho _{1}+\cdots +\varrho _{d}$\emph{,} for distinct
rays $\varrho _{1},\ldots ,\varrho _{d}$. The \textit{multiplicity} mult$%
\left( \sigma ;N\right) $ of $\sigma $ with\emph{\ r}espect to $N$ is
defined as
\begin{equation*}
\text{mult}\left( \sigma ;N\right) :=\frac{\text{det}\left( \mathbb{Z}%
\,n\left( \varrho _{1}\right) \oplus \cdots \oplus \mathbb{Z}\,n\left(
\varrho _{d}\right) \right) }{\text{det}\left( N_{\sigma }\right) },
\end{equation*}
where $N_{\sigma }$ is the lattice in lin$\left( \sigma \right) $ induced by
$N.$ If\emph{\ }mult$\left( \sigma ;N\right) =1$\emph{, }then\emph{\ }$%
\sigma $\emph{\ }is called a \textit{basic cone} w.r.t. $N$.
\end{definition}

\begin{theorem}[Smoothness criterion]
\label{SMCR}The affine toric variety $U_{\sigma }$ is smooth iff $\sigma $
is basic \textit{w.r.t.} $N$. \emph{(}Correspondingly, an arbitrary toric
variety $X\left( N,\Delta \right) $ is smooth if and only if it is
simplicial and each s.c.p. cone $\sigma \in \Delta $ is basic \textit{w.r.t.}
$N$.\emph{)}
\end{theorem}

\noindent \textit{Proof. }See \cite[Ch. I, Thm. 4, p. 14]{KKMS} and
\cite[Thm. 1.10, p. 15]{Oda}.\hfill $\square $

\begin{theorem}[$\mathbb{Q}$-factoriality]
A toric variety $X\left( N,\Delta \right) $ is $\mathbb{Q}$-factorial if and
only if $\Delta $ is simplicial, i.e., if and only if $X\left( N,\Delta
\right) $ has at most abelian quotient singularities.
\end{theorem}

\noindent {}\textit{Proof}. Since this is a local property, it is enough to
consider the case in which $X\left( N,\Delta \right) =U_{\sigma },$ where
the cone $\sigma =\mathbb{R}_{\geq 0}v_{1}+\cdots +\mathbb{R}_{\geq 0}v_{r}$
is of maximal dimension, Gen$\left( \sigma \right) =\{v_{1},...,v_{r}\},$
and $\Delta =\left\{ \tau \ \left| \ \tau \preceq \sigma \right. \right\} $.
$U_{\sigma }$ is $\mathbb{Q}$-factorial if and only if all the $T_{N}$%
-invariant prime divisors $D_{v_{i}}$ are $\mathbb{Q}$-Cartier. This is
equivalent to the existence of elements $m_{i}\in M_{\mathbb{Q}}=M\otimes _{%
\mathbb{Z}}\mathbb{Q},$ with $M=$ Hom$_{\mathbb{Z}}\left( N,\mathbb{Z}%
\right) ,$ for which $\left\langle m_{i},v_{j}\right\rangle =\delta _{ij}$
(the Kronecker delta). But this means that $\sigma $ is simplicial.\hfill $%
\square \bigskip $\newline
Next theorem is due to Stanley (\cite[\S 6]{Stanley1}), who worked directly
with the monoidal $\mathbb{C}$-algebra $\mathbb{C}\left[ M\cap \sigma ^{\vee
}\right] $, as well as to Ishida (\cite[\S 7]{Ishida}), Danilov and Reid (
\cite[p. 294]{REID1}), who provided a purely algebraic-geometric
characterization of the Gorenstein property.

\begin{theorem}[Gorenstein property]
\label{GOR-PR}Let $N$ be a free $\mathbb{Z}$-module of rank $r$, and $\sigma
\subset N_{\mathbb{R}}$ \ a s.c.p. cone of dimension $d\leq r$. Then the
following conditions are equivalent\emph{:\medskip }\newline
\emph{(i)} $\ \ U_{\sigma }$ is Gorenstein.\medskip\ \newline
\emph{(ii)} \ There exists an element $m_{\sigma }$ of $M$, such that
\begin{equation*}
M\cap \left( \text{\emph{int}}\left( \sigma ^{\vee }\right) \right)
=m_{\sigma }+M\cap \sigma ^{\vee }.\medskip \newline
\end{equation*}
$\smallskip $\emph{(iii)} \emph{Gen}$\left( \sigma \right) \subset \mathbf{H}
$, where $\mathbf{H}$ denotes a primitive affine hyperplane of \ $\left(
N_{\sigma }\right) _{\mathbb{R}}$ .\medskip \newline
Moreover, if \emph{\ }$d=r$, then $m_{\sigma }$ in \emph{(ii) }is a \emph{%
uniquely determined} primitive element of $M\cap \left( \text{\emph{int}}%
\left( \sigma ^{\vee }\right) \right) $ and $\mathbf{H}$ in \emph{(iii) }%
equals
\begin{equation*}
\mathbf{H}=\left\{ \mathbf{x}\in N_{\mathbb{R}}\ \left| \ \left\langle
m_{\sigma },\mathbf{x}\right\rangle =1\right. \right\} .
\end{equation*}
\end{theorem}

\begin{definition}
\label{LP}\emph{\ }If $N_{1}$ and $N_{2}$\ are two free $\mathbb{Z}$-modules
(not necessarily of\emph{\ }the same rank) and $P_{1}\subset \left(
N_{1}\right) _{\mathbb{R}}$,\emph{\ }$P_{2}\subset \left( N_{2}\right) _{%
\mathbb{R}}$\emph{\ }two lattice polytopes w.r.t. them, we shall say that $%
P_{1}$\ and $P_{2}$\emph{\ }are \textit{lattice equivalent }to each other,
if $P_{1}$\ is affinely equivalent to $P_{2}$\ via an\emph{\ }affine map%
\emph{\ }$\varpi :\left( N_{1}\right) _{\mathbb{R}}\rightarrow \left(
N_{2}\right) _{\mathbb{R}}$\emph{, }such that $\varpi \left| _{\text{aff}%
\left( P\right) }\right. :$ aff$\left( P\right) \rightarrow $ aff$\left(
P^{\prime }\newline
\right) $ is a bijection mapping\emph{\ }$P_{1}$\emph{\ }onto the
(necessarily equidimensional) polytope\emph{\ }$P_{2}$, every $i$%
-dimensional face of\emph{\ }$P_{1}$\emph{\ }onto an\emph{\ }$i$-dimensional%
\emph{\ }face of\emph{\ }$P_{2}$, for all\emph{\ }$i$, $0\leq i\leq $ dim$%
\left( P_{1}\right) =$ dim$\left( P_{2}\right) $, and, in addition, $%
N_{P_{1}}$\emph{\ }onto the lattice\emph{\ }$N_{P_{2}}$, where by\emph{\ }$%
N_{P_{j}}$\emph{\ }is\emph{\ }meant the sublattice of\emph{\ }$N_{j}$
generated (as subgroup) by aff$\left( P_{j}\right) \cap N_{j}$\emph{, }$%
j=1,2.$\emph{\ }
\end{definition}

\begin{definition}
\label{NAKAJI}\textit{Nakajima polytopes} $P\subset \mathbb{R}^{r}$ are
lattice polytopes w.r.t. the usual rectangular lattice $\mathbb{Z}^{r}$,
defined inductively as follows: In dimension $0,$ $P$ is an (arbitrary)
point of $\mathbb{R}^{r},$ while in dimension $d\leq r,$ $P$ is of the form
\begin{equation*}
P=\left\{ \left. \mathbf{x}=\left( \mathbf{x}^{\prime },x_{r}\right) \in
F\times \mathbb{R\subset \mathbb{R}}^{r}\mathbb{\,}\right| \,0\leq x_{r}\leq
\left\langle \mathbf{m},\mathbf{x}^{\prime }\right\rangle \right\} ,
\end{equation*}
where the facet $F\subset \mathbb{R}^{r-1}$ is a Nakajima polytope of
dimension $d-1$, and $\mathbf{m\in }\left( \mathbb{Z}^{r}\right) ^{\vee }$
is a functional taking non-negative values on $F.$
\end{definition}

\begin{theorem}[Toric L.C.I.'s]
Let $N$ be a free $\mathbb{Z}$-module of rank $r$, and $\sigma \subset N_{%
\mathbb{R}}$ \ a s.c.p. cone of dimension $d\leq r$, such that $U_{\sigma }$
is Gorenstein. Writing $\sigma $ as $\sigma =\sigma ^{\prime }\oplus \left\{
\mathbf{0}\right\} $ with $\sigma ^{\prime }$ a $d$-dimensional cone in $%
\left( N_{\sigma }\right) _{\mathbb{R}}$, we obtain an analytic isomorphism:
\begin{equation*}
U_{\sigma }\cong U_{\sigma ^{\prime }}\times \left( \mathbb{C}^{\ast
}\right) ^{r-d}.
\end{equation*}
Let $m_{\sigma ^{\prime }}$ be the unique primitive element of $M_{\sigma
}\cap \left( \text{\emph{int}}\left( \left( \sigma ^{\prime }\right) ^{\vee
}\right) \right) $, as it is defined in \emph{Theorem \ref{GOR-PR}}. Then $%
U_{\sigma }$ is a \emph{local complete intersection\footnote{%
Obviously, for $d=r$, $U_{\sigma }$ is a ``g.c.i.'' in the sense of \S \ref
{SECTION1}.}\ }if and only if the lattice polytope \textbf{\ }\
\begin{equation*}
P_{\sigma ^{\prime }}:=\sigma ^{\prime }\cap \left\{ \mathbf{x}\in \left(
N_{\sigma }\right) _{\mathbb{R}}\ \left| \ \right. \left\langle m_{\sigma
^{\prime }},\mathbf{x}\right\rangle =1\right\}
\end{equation*}
is lattice equivalent to a Nakajima polytope \emph{(}cf.\emph{\ \ref{LP},
\ref{NAKAJI})}
\end{theorem}

\noindent {}{}This was proved in \cite{Ishida} for dimension 3 and in \cite
{Nakajima} for arbitrary dimensions.\bigskip

\noindent \textsf{(i)}\textit{\ }A \textit{map of fans\ }$\varpi :\left(
N^{\prime },\Delta ^{\prime }\right) \rightarrow \left( N,\Delta \right) $
is a $\mathbb{Z}$-linear homomorphism $\varpi :N^{\prime }\rightarrow N$
whose scalar extension $\varpi \otimes _{\mathbb{Z}}$id$_{\mathbb{R}}:N_{%
\mathbb{R}}^{\prime }\rightarrow N_{\mathbb{R}}$ satisfies the property:
\begin{equation*}
\forall \sigma ^{\prime },\ \sigma ^{\prime }\in \Delta ^{\prime }\ \ \text{
}\exists \ \sigma ,\ \sigma \in \Delta \ \ \text{ with\ \ }\varpi \left(
\sigma ^{\prime }\right) \subset \sigma \,.
\end{equation*}
$\varpi \otimes _{\mathbb{Z}}$id$_{\mathbb{C}^{\ast }}:T_{N^{\prime
}}=N^{\prime }\otimes _{\mathbb{Z}}\mathbb{C}^{\ast }\rightarrow
T_{N}=N\otimes _{\mathbb{Z}}\mathbb{C}^{\ast }$ is a homomorphism from $%
T_{N^{\prime }}$ to $T_{N}$ and the scalar extension $\varpi ^{\vee }\otimes
_{\mathbb{Z}}$id$_{\mathbb{R}}:M_{\mathbb{R}}\rightarrow M_{\mathbb{R}%
}^{\prime }$ of the dual $\mathbb{Z}$-linear map $\varpi ^{\vee
}:M\rightarrow M^{\prime }$ induces canonically an \textit{equivariant
holomorphic map }
\begin{equation*}
\varpi _{\ast }:X\left( N^{\prime },\Delta ^{\prime }\right) \rightarrow
X\left( N,\Delta \right) \,.
\end{equation*}
This map is\textit{\ proper} if and only if $\varpi ^{-1}\left( \left|
\Delta \right| \right) =\left| \Delta ^{\prime }\right| .$ In particular, if
$N=N^{\prime }$ and $\Delta ^{\prime }$ is a refinement of $\Delta $, then id%
$_{\ast }:X\left( N,\Delta ^{\prime }\right) \rightarrow X\left( N,\Delta
\right) $ is \textit{proper}\emph{\ }and \textit{birational }(cf. \cite[Thm.
1.15 and Cor. 1.18]{Oda}).\bigskip

\noindent \textsf{(j)}\textit{\ }By Carath\'{e}odory's Theorem concerning
convex polyhedral cones (cf. \cite[III 2.6 and V 4.2]{Ewald}) one can choose
a refinement $\Delta ^{\prime }$ of any given fan $\Delta $, so that $\Delta
^{\prime }$ becomes simplicial. Since further subdivisions of $\Delta
^{\prime }$ reduce the multiplicities of its cones, we may arrive (after
finitely many subdivisions) at a fan $\widetilde{\Delta }$ having only basic
cones. Hence, for every toric variety $X\left( N,\Delta \right) $ there
exists a refinement $\widetilde{\Delta }$ of $\Delta $ consisting of
exclusively basic cones w.r.t. $N$, i.e., such that
\begin{equation*}
f=\text{id}_{\ast }:X(N,\widetilde{\Delta })\longrightarrow X\left( N,\Delta
\right)
\end{equation*}
is a $T_{N}$-equivariant (full) desingularization.$\medskip $

\begin{theorem}
\label{LOG-TERMIN}All $\mathbb{Q}$-Gorenstein toric varieties have at most
log-terminal singularities.
\end{theorem}

\noindent {}\textit{Proof}. We may again assume that $X=X\left( N,\Delta
\right) =U_{\sigma }$ (where $\sigma =\mathbb{R}_{\geq 0}v_{1}+\cdots +%
\mathbb{R}_{\geq 0}v_{r}$ is a cone of maximal dimension, with Gen$\left(
\sigma \right) =\{v_{1},...,v_{r}\}$ and $\Delta =\left\{ \tau \ \left| \
\tau \preceq \sigma \right. \right\} $). Since $\omega _{X}=\mathcal{O}%
_{X}\left( K_{X}\right) $ with $K_{X}=-\sum D_{v_{i}}\leq 0,$ for $K_{X}$ to
be $\mathbb{Q}$-Cartier means that $K_{X}$ is a ($T_{N}$-invariant) Cartier
divisor after multiplication by an integer. This multiple of $K_{X}$ has to
be a divisor of the form div$\left( \mathbf{e}(u)\right) ,$ for some $u\in M$
(cf. \cite[Prop. 2.1, pp. 68-69]{Oda} or \cite[Lemma of p. 61]{Fulton}).
Hence, there must be an $m_{\sigma }\in M_{\mathbb{Q}},$ such that $%
\left\langle m_{\sigma },\sigma \right\rangle \geq 0$ and $\left\langle
m_{\sigma },v_{j}\right\rangle =1$ (with $m_{\sigma }$ regarded as a \textit{%
linear} support function on $\sigma $, cf. \textsf{(g)}). Let now
\begin{equation*}
f:Y=X\left( N,\Delta ^{\prime }\right) \rightarrow X=X(N,\Delta )=U_{\sigma }
\end{equation*}
be a desingularization of $X$ obtained by a subdivision of $\sigma $ into
smaller basic strongly convex rational polyhedral cones. Suppose that the
primitive lattice points of the \textit{new} introduced rays in $\Delta
^{\prime }$ are $v_{1}^{\prime },v_{2}^{\prime },\ldots ,v_{s}^{\prime }.$
Then
\begin{equation*}
f^{\ast }K_{X}=-\sum_{v\in \text{Gen}\left( \Delta ^{\prime }\right)
}\left\langle m_{\sigma },v\right\rangle D_{v}.
\end{equation*}
Since
\begin{equation*}
\mathcal{O}_{Y}\left( K_{Y}\right) \left| _{\text{Reg}\left( Y\right)
}\right. =\mathcal{O}_{Y}\left( f^{\ast }K_{X}\right) \left| _{\text{Reg}%
\left( Y\right) }\right.
\end{equation*}
we have
\begin{equation}
K_{Y}-f^{\ast }K_{X}=\sum_{j=1}^{s}\left( \left\langle m_{\sigma
},v_{j}^{\prime }\right\rangle -1\right) D_{v_{j}^{\prime }}
\label{DISCREPANCY}
\end{equation}
Since the discrepancy is $>-1,$ $\left( U_{\sigma },\text{orb}\left( \sigma
\right) \right) $ is a log-terminal singularity whose index equals $\ell =$
min$\{\kappa \in \mathbb{Z}_{\geq 1}\ \left| \ m_{\sigma }\kappa \in \right.
M\}.$\hfill $\square $

\begin{theorem}
\emph{\ }All toric singularities are rational singularities.
\end{theorem}

\noindent {}\textit{Proof. }\ For $\mathbb{Q}$-Gorenstein toric
singularities this follows from Theorems \ref{LOG-RAT} and \ref{LOG-TERMIN}.
For the general case see Ishida \cite{Ishida} and Oda \cite[Cor. 3.9, p. 125]
{Oda}.\hfill $\square \medskip $

\noindent {}Formula (\ref{DISCREPANCY}) gives us the following purely
combinatorial characterization of terminal (resp. canonical) toric
singularities.

\begin{corollary}
A toric singularity $\left( U_{\sigma },\text{\emph{orb}}\left( \sigma
\right) \right) $ is terminal \emph{(}resp. canonical\emph{)} of index $\ell
$ if and only if{\small
\begin{equation*}
\begin{array}{c}
\exists \ m_{\sigma }\in M:\left\langle m_{\sigma },\text{\emph{Gen}}\left(
\sigma \right) \right\rangle =\ell \text{ }\emph{and\ }\left\langle
m_{\sigma },n\right\rangle >\ell \text{\emph{\ \ for all} \ }n\in \sigma
\cap N\mathbb{r}(\{\mathbf{0}\}\cup \text{\emph{Gen}}\left( \sigma \right) ),
\\
\  \\
\text{\emph{(}resp. }\exists \ m_{\sigma }\in M:\left\langle m_{\sigma },%
\text{\emph{Gen}}\left( \sigma \right) \right\rangle =\ell \text{ }\emph{%
and\ \ }\left\langle m_{\sigma },n\right\rangle \geq \ell \,\text{\emph{\
for all} \ }n\in \sigma \cap N\mathbb{r}\{\mathbf{0}\}\text{\emph{).}}
\end{array}
\end{equation*}
}\emph{(}$m_{\sigma \text{ }}$\emph{\ is uniquely determined whenever }$%
\sigma $\emph{\ is of maximal dimension in the fan).}
\end{corollary}

\newpage \noindent \textsf{(k) }\textbf{Recapitulation}. This is divided
into two parts. The first one contains the standard dictionary: \medskip
{\footnotesize
\begin{equation*}
\begin{tabular}{|l|l|}
\hline
\textbf{Discrete Geometry} & \textbf{Algebraic Geometry} \\ \hline\hline
$
\begin{array}{l}
\text{strongly convex rational} \\
\text{pol.cone }\sigma \subset N_{\mathbb{R}}\cong \mathbb{R}^{r}
\end{array}
$ & $
\begin{array}{l}
\text{affine toric variety} \\
U_{\sigma }=\text{Spec}(\mathbb{C}[\sigma ^{\vee }\cap M])
\end{array}
$ \\ \hline
Fans $\Delta $ & $
\begin{array}{l}
\text{toric varieties }X\left( N,\Delta \right) \\
\text{(after glueing)}
\end{array}
$ \\ \hline
$\tau $ face of $\sigma $ & $U_{\tau }$ open subset of $U_{\sigma }$ \\
\hline
$\left| \Delta \right| =\mathbb{R}^{r}$ & $X\left( N,\Delta \right) $ is
complete \\ \hline
All cones of $\Delta $ are basic & $X\left( N,\Delta \right) $ is
non-singular \\ \hline
All cones of $\Delta $ are simplicial & $
\begin{array}{l}
X\left( N,\Delta \right) \text{ has at most abelian} \\
\text{quotient singularities}
\end{array}
$ \\ \hline
$
\begin{array}{l}
\Delta \text{ admits a strictly convex} \\
\text{upper support function}
\end{array}
$ & $
\begin{array}{l}
X\left( N,\Delta \right) \text{ is a quasi-} \\
\text{projective variety}
\end{array}
$ \\ \hline
$\Delta ^{\prime }$ a cone subdivision of $\Delta $ & $
\begin{array}{l}
X\left( N,\Delta ^{\prime }\right) \longrightarrow X\left( N,\Delta \right)
\text{ proper} \\
\text{birational morphism}
\end{array}
$ \\ \hline
\end{tabular}
\bigskip
\end{equation*}
} \bigskip \noindent {}while the second one describes the main properties of
toric singularities: \medskip {\footnotesize
\begin{equation*}
\begin{tabular}{|l|l|}
\hline
\textbf{Discrete Geometry} & \textbf{Algebraic Geometry} \\ \hline\hline
$
\begin{array}{l}
\text{strongly convex rational} \\
\text{pol.cone }\sigma \subset \mathbb{R}^{r}
\end{array}
$ & $
\begin{array}{l}
\, \\
U_{\sigma }\text{ is always normal}, \\
\text{and Cohen-Macaulay} \\
\,
\end{array}
$ \\ \hline
$
\begin{array}{c}
\  \\
\begin{array}{l}
\text{orb}\left( \sigma \right) \in U_{\sigma },\ \sigma \text{ non-basic
cone} \\
\text{of maximal dimension }r \\
\text{(keep this assumption in what follows)}
\end{array}
\\
\
\end{array}
$ & $
\begin{array}{l}
\left( U_{\sigma },\text{orb}\left( \sigma \right) \right) \text{ is a} \\
\text{rational singularity}
\end{array}
$ \\ \hline
$
\begin{array}{l}
\exists !\text{ }m_{\sigma }\in M_{\mathbb{Q}}: \\
\text{Gen}\left( \sigma \right) \subset \left\{ y\in N_{\mathbb{R}}\ \left|
\ \right. \left\langle m_{\sigma },y\right\rangle =1\right\}
\end{array}
$ & $
\begin{array}{l}
\, \\
U_{\sigma }\text{ is }\mathbb{Q}\text{-Gorenstein} \\
\text{and i.p. }\left( U_{\sigma },\text{orb}\left( \sigma \right) \right)
\text{ is a } \\
\text{log-terminal singularity of index} \\
\text{min}\{\kappa \in \mathbb{Z}_{\geq 1}\ \left| \ m_{\sigma }\kappa \in
\right. M\} \\
\,
\end{array}
$ \\ \hline
$\left\{
\begin{array}{l}
\exists !\text{ }m_{\sigma }\in M_{\mathbb{Q}}: \\
\text{Gen}\left( \sigma \right) \subset \left\{ y\in N_{\mathbb{R}}\ \left|
\ \right. \left\langle m_{\sigma },y\right\rangle =1\right\} \\
\text{and \ \ }N\cap \sigma \cap \left\{ y\in N_{\mathbb{R}}\ \left| \
\right. \left\langle m_{\sigma },y\right\rangle <1\right\} =\{\mathbf{0}\}\
\end{array}
\right. $ & $
\begin{array}{l}
\, \\
\left( U_{\sigma },\text{orb}\left( \sigma \right) \right) \text{ is a} \\
\text{canonical singularity of index} \\
\text{min}\{\kappa \in \mathbb{Z}_{\geq 1}\ \left| \ m_{\sigma }\kappa \in
\right. M\} \\
\,
\end{array}
$ \\ \hline
$\left\{
\begin{array}{l}
\exists !\text{ }m_{\sigma }\in M_{\mathbb{Q}}: \\
\text{Gen}\left( \sigma \right) \subset \left\{ \mathbf{x}\in N_{\mathbb{R}%
}\ \left| \ \right. \left\langle m_{\sigma },\mathbf{x}\right\rangle
=1\right\} \\
\text{and \ \ }N\cap \sigma \cap \left\{ \mathbf{x}\in N_{\mathbb{R}}\
\left| \ \right. \left\langle m_{\sigma },\mathbf{x}\right\rangle \leq
1\right\} \\
=\{\mathbf{0}\}\cup \text{Gen}\left( \sigma \right) \
\end{array}
\right. $ & $
\begin{array}{l}
\, \\
\left( U_{\sigma },\text{orb}\left( \sigma \right) \right) \text{ is a} \\
\text{terminal singularity of index} \\
\text{min}\{\kappa \in \mathbb{Z}_{\geq 1}\ \left| \ m_{\sigma }\kappa \in
\right. M\} \\
\,
\end{array}
$ \\ \hline
$
\begin{array}{c}
\  \\
\left\{
\begin{array}{l}
\exists !\text{ }m_{\sigma }\in M: \\
\text{Gen}\left( \sigma \right) \subset \left\{ \mathbf{x}\in N_{\mathbb{R}%
}\ \left| \ \right. \left\langle m_{\sigma },\mathbf{x}\right\rangle
=1\right\} \\
\ \text{with }\sigma \text{ supporting the lattice polytope} \\
P_{\sigma }:=\sigma \cap \left\{ \mathbf{x}\in N_{\mathbb{R}}\ \left| \
\right. \left\langle m_{\sigma },\mathbf{x}\right\rangle =1\right\} \text{
and} \\
N\cap \text{conv}\left( \{\mathbf{0}\}\ \cup P_{\sigma }\right) =\{\mathbf{0}%
\}\ \cup \left( N\cap P_{\sigma }\right)
\end{array}
\right. \\
\
\end{array}
$ & $
\begin{array}{l}
\left( U_{\sigma },\text{orb}\left( \sigma \right) \right) \text{ is a } \\
\text{Gorenstein singularity} \\
\text{(i.e., canonical of \ index \ }1\text{)}
\end{array}
$ \\ \hline
$
\begin{array}{l}
P_{\sigma }\text{ (as above) is lattice-equivalent} \\
\text{to a Nakajima polytope}
\end{array}
$ & $
\begin{array}{l}
\, \\
\left( U_{\sigma },\text{orb}\left( \sigma \right) \right) \text{ is a
locally} \\
\text{complete intersection } \\
\text{singularity} \\
\,
\end{array}
$ \\ \hline
\end{tabular}
\end{equation*}
} 

\begin{remark}
(i) A lattice polytope $P$\ is called \textit{elementary} if the lattice
points belonging to it are exactly its vertices. A lattice simplex is said
to be \textit{basic} (or \textit{unimodular}) if its vertices constitute a
part of a\emph{\ }$\mathbb{Z}$-basis of the reference lattice (or
equivalently, if its relative, normalized volume equals $1$).\smallskip
\newline
(ii) It is clear by the fifth row of the above table that the classification
of terminal (resp. $\mathbb{Q}$-factorial terminal) $r$-dimensional
Gorenstein toric singularities is \textit{equivalent} to the classification
of all $(r-1)$-dimensional elementary lattice polytopes (resp. elementary
lattice simplices). The readers, who would like to learn more about partial
classification results (on both ``sides'') in dimensions $\geq 3$, are
referred (in chronological order) to \cite{DANILOV}, \cite{Mo-Ste}, \cite
{MMM}, \cite{ISHIDA-IWASHITA}, \cite[pp. 34-36]{Oda}, \cite{Sankaran}, \cite
{BORISOV1}, and \cite{BORISOV2}, for investigations from the point of view
of computational algebraic geometry, as well as to \cite{KANTOR}, \cite{B-K}%
, \cite{Haase-Ziegler}, for a study of ``width-functions'' (which are
closely related to the index of the corresponding singularities) from the
point of view of number theory, combinatorics and integer programming%
\footnote{%
Note that in the literature, instead of ``elementary'' polytopes, there are
also in use \textit{different names} like ``fundamental'',
``lattice-point-free'', ``hollow'', or even ``empty'' polytopes.}.\smallskip
\newline
(iii) On the other extreme, all toric l.c.i.-singularities admit \textit{%
crepant} resolutions in \textit{all }dimensions. (This was proved recently
in \cite{DHaZ} by showing inductively that all Nakajima polytopes admit
lattice triangulations consisting exclusively of basic simplices; cf.
Theorem \ref{SMCR}). Nevetheless, regarding non-l.c.i., Gorenstein toric
singularities, the determination of necessary and sufficient ``intrinsic''
conditions, under which they possess resolutions of this kind, remains an
unsolved problem. Finally, let us stress that one can always reduce
log-terminal non-canonical toric singularities to canonical ones by a
torus-equivariant, natural ``canonical modification'' which is uniquely
determined (cf. below footnote to \ref{MAIN} (i)).
\end{remark}

\section{Toric two- and three-dimensional singularities, \newline
and their resolutions\label{SECTION4}}

\noindent {}Now we focus on the desingularization methods of low-dimensional
toric singularities.\medskip\

\noindent {}$\bullet $ All \textbf{2-dimensional} toric singularities $%
\left( U_{\sigma },\text{orb}\left( \sigma \right) \right) $ are abelian (in
fact, cyclic) quotient singularities and can be treated by means of the
finite continued fractions (see \cite[\S 1.6]{Oda}). For them there exists
always a uniquely determined ``minimal'' resolution\footnote{%
This is actually a ``good minimal'' resolution in the sense of \S \ref
{SECTION2}.}. In fact, this resolution has the nice property that
\begin{equation*}
\mathbf{Hilb}_{N}\left( \sigma \right) =\text{Gen}\left( \Delta ^{\prime
}\right) ,
\end{equation*}
where $\Delta ^{\prime }$ is the fan which refines $\sigma $ into basic
cones and $\mathbf{Hilb}_{N}\left( \sigma \right) $ the Hilbert basis of $%
\sigma $ w.r.t. $N.$ Figure \textbf{1} shows this minimal resolution $%
U_{\sigma }\overset{f}{\longleftarrow }X\left( N,\Delta ^{\prime }\right) $
of $U_{\sigma }$ for the cone \
\begin{equation*}
\sigma =\text{pos}\left( \left\{ \left( 1,0\right) ,(4,5)\right\} \right)
\subset \mathbb{R}^{2}
\end{equation*}
(w.r.t. the standard lattice $N=\mathbb{Z}^{2}$) constructed by its
subdivision into two basic subcones.

\begin{figure}[h]
\begin{center}
\includegraphics[width=2.5in,height=2.5in]{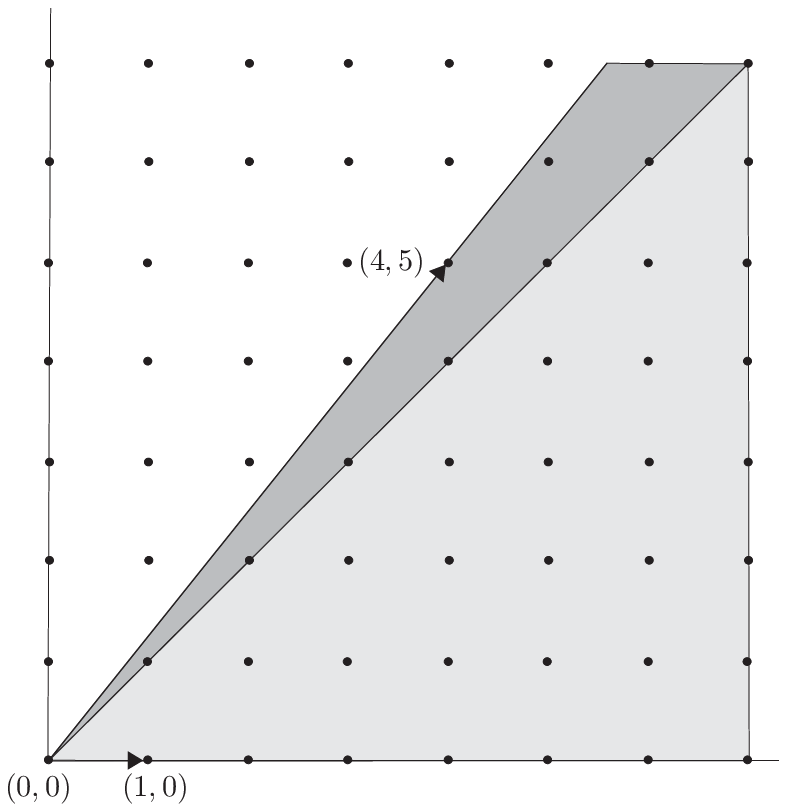} \vspace{0.3cm} \\[0pt]
\textbf{Fig. 1}
\end{center}
\end{figure}

\noindent Note that $U_{\sigma }=\mathbb{C}^{2}/G$, where $G\subset $ GL$(2,%
\mathbb{C})$ is the finite cyclic group generated by the matrix
\begin{equation*}
\left(
\begin{array}{cc}
e^{\frac{2\pi \sqrt{-1}}{5}} & 0 \\
0 & e^{\frac{2\pi \sqrt{-1}}{5}}
\end{array}
\right) \ ,
\end{equation*}
and
\begin{equation*}
\sigma ^{\vee }=\text{pos}\left( \left\{ \left( 0,1\right) ,\left(
5,-4\right) \right\} \right) \subset \left( \mathbb{R}^{2}\right) ^{\vee }.
\end{equation*}
In particular, $\mathbb{C}[\left( \mathbb{R}^{2}\right) ^{\vee }\cap \sigma
^{\vee }]\cong \mathbb{C}\left[ u,w\right] ^{G}$ is generated by the
monomials
\begin{equation*}
\{\left. u^{i}w^{j}\,\right| \,i+j \equiv 0(\text{mod} \mathrm{5})\}.
\end{equation*}
The \textit{minimal} generating system of $(\mathbb{R}^{2})^{\vee }\cap
\sigma ^{\vee }$ equals
\begin{equation*}
\mathbf{Hilb}_{(\mathbb{R}^{2})^{\vee }}\left( \sigma ^{\vee }\right)
=\{k_{0},k_{1},k_{2},k_{3},k_{4},k_{5}\},
\end{equation*}
with $k_{0}=\left( 0,1\right) $, $k_{5}=\left( 5,-4\right) $ the primitive
vectors spanning $\sigma ^{\vee }.$ The remaining elements are determined by
the vectorial matrix multiplication
\begin{equation*}
\left(
\begin{array}{c}
k_{1} \\
k_{2} \\
k_{3} \\
k_{4}
\end{array}
\right) =\left(
\begin{array}{cccc}
2 & -1 & 0 & 0 \\
-1 & 2 & -1 & 0 \\
0 & -1 & 2 & -1 \\
0 & 0 & -1 & 2
\end{array}
\right) ^{-1}\left(
\begin{array}{c}
k_{0} \\
(0,0) \\
(0,0) \\
k_{5}
\end{array}
\right)
\end{equation*}
where the diagonal entries of the first matrix of the right-hand side are
exactly those arising in the continued fraction development
\begin{equation*}
\frac{5}{4}=2-\frac{1}{2-\frac{1}{2-\frac{1}{2}}}\ .
\end{equation*}
(cf. \cite[Prop. 1.21, p. 26]{Oda}). $U_{\sigma }$ has embedding dimension $%
6 $, and taking into account the linear dependencies between the members of
this Hilbert basis, it can be described as the zero-locus of $\binom{5}{2}%
=\allowbreak 10$ square-free binomials (by applying Proposition \ref{EMB}).
More precisely (cf. \cite{RIEMENSCHNEIDER}), setting $z_{i}:=\mathbf{e}%
\left( k_{i}\right) $, $0\leq i\leq 5$, we obtain:\smallskip\
\begin{equation*}
U_{\sigma }\cong \left\{ \left( z_{0},z_{1},z_{2},z_{3},z_{4},z_{5}\right)
\in \mathbb{C}^{6}\,\left| \,\text{rank}\left(
\begin{array}{ccccc}
z_{0} & z_{1} & z_{2} & z_{3} & z_{4} \\
z_{1} & z_{2} & z_{3} & z_{4} & z_{5}
\end{array}
\right) \right. \leq 1\right\} .
\end{equation*}
\bigskip

\noindent {}$\bullet $ In \textbf{dimension 3} things are more complicated.
The singularities can be resolved by more or less ``canonical'' procedures
(and by \textit{projective} birational morphisms) but the ``uniqueness'' is
mostly lost (even if one requires ``minimality'' in the sense of ``MMP'',
i.e., $\mathbb{Q}$-terminalizations), though ``not completely'' (of course,
as usual, meant ``up to isomorphisms in codimension 1''). In the literature
you may find several different approaches:\medskip

\noindent {}\textbf{(a)} \textbf{Hilb-desingularizations.} Bouvier and
Gonzalez-Sprinberg used in \cite{B-G-S1}, \cite{B-G-S2} ``$\mathbf{Hilb}$%
\textbf{-}desingularizations'' (i.e., again with $\mathbf{Hilb}_{N}\left(
\sigma \right) =$ Gen$\left( \Delta ^{\prime }\right) $) to resolve
three-dimensional toric singularities $\left( U_{\sigma },\text{orb}\left(
\sigma \right) \right) .$ Their exposition is well-structured and fits
together with more general concepts of MMP (see \cite{DANILOV}, \cite
{KAWAMATA}, \cite{MORI2}). Nevertheless, any kind of ``uniqueness'' is (in
general) lost already from the ``second'' step by producing their ``minimal
terminal subdivisions''.\medskip

\noindent {}\textbf{(b) Hilb-desingularizations with extra lexicographic
ordering. }Similar constructive method, due to Aguzzoli and Mundici\textsc{\
}\cite{AG-MU} (again by making use of ``indispensable'' exceptional
divisors), but with a fixed ordering for performing the starring operations.
The ``uniqueness'' is lost as long as one uses other orderings. Moreover,
they also assume $\mathbb{Q}$-factorialization from the very
beginning.\medskip

\noindent {}\textbf{(c) Distinguished crepant desingularization (via the
Hilbert scheme of }$G$\textbf{-clusters) for the case in which }$\sigma $%
\textbf{\ is simplicial and }$U_{\sigma }=\mathbb{C}^{3}/G$ \ \textbf{%
Gorenstein}. Motivated by the so-called ``McKay correspondence'' in
dimension 3 (see \cite{REID5} for a recent exposition), Nakamura \cite
{Nakamura} and Craw \& Reid \cite{CRAW1}, \cite{CRAW2}, \cite{Craw-Reid}
have constructed a \textit{distinguished} crepant resolution of $U_{\sigma }=%
\mathbb{C}^{3}/G$ \ ($G$ an abelian finite subgroup of SL($3,\mathbb{C}$)),
by expressing the $G$-orbit Hilbert scheme of $\mathbb{C}^{3}$ as a fan of
basic cones supporting certain ``regular tesselations'' of the corresponding
junior simplex. Unfortunately, it is not known if there is an analogue of
this method for the case in which $\sigma $ is not simplicial, and it is
known that it does not work (at least as crepant resolution) in higher
dimensions.\medskip

\noindent {}\textbf{(d) Desingularizations via the ``initial strategy''of M.
Reid }\cite{REID1}, \cite{REID2}, \cite{REID4}\textbf{. }This method was
already described in broad outline in section \ref{SECTION2}. Applying it
for arbitrary three-dimensional toric singularities (cf. \cite{DHZ}), we win
more explicit information about the resolution (compared with the general
case\footnote{%
Cf. \cite[Rem. 6.10, p. 308]{REID1} and \cite[Rem. 0.8 (d), p. 135, and Ex.
2.7, p. 145]{REID2}.\smallskip}) by obtaining the following:

\begin{theorem}
\label{MAIN}Let $\left( U_{\sigma },\text{\emph{orb}}\left( \sigma \right)
\right) $ be an arbitrary 3-dimensional toric singularity, where $U_{\sigma
}=X\left( N,\Delta \right) $. \ \medskip \newline
\emph{(i)} There exists a \emph{uniquely determined} partial
desingularization
\begin{equation}
U_{\sigma }=X\left( N,\Delta \right) \longleftarrow X\left( N,\Delta _{\text{%
\emph{can}}}\right)  \label{CAN-MOD}
\end{equation}
such that $X\left( N,\Delta _{\text{\emph{can}}}\right) $ has only canonical
singularities. \emph{(}This is actually a general fact valid for all
dimensions\emph{\footnote{%
To construct the ``canonical modification'' (\ref{CAN-MOD}) one has just to
pass to the fan-subdivision
\begin{equation*}
\Delta _{\text{can}}=\left\{ \left\{ \left. \text{pos}(F)\right| \text{ }F%
\text{ faces of conv}(\sigma \cap N\mathbb{r}\{\mathbf{0}\})\right\} ,\
\sigma \in \Delta \right\}
\end{equation*}
}).\smallskip }\newline
\emph{(ii)} The singularities of $X\left( N,\Delta _{\text{\emph{can}}%
}\right) $ which have index $>1$ can be reduced to canonical singularities
of index $1$ up to cyclic coverings \emph{(}via lattice dilations\emph{).}
Moreover, one may treat them in a special manner \emph{(}either by the
Ishida-Iwashita classification \emph{\cite{ISHIDA-IWASHITA}} or by Bouvier,
Gonzalez-Sprinberg's \emph{\cite{B-G-S2}} $\mathbf{Hilb}$\textbf{-}%
desingularizations\emph{\footnote{%
For the ``uniqueness'' of this reduction procedure to the index $1$ case
(from the point of view of toric geometry) see, in particular,
\cite[Th\'{e}or\`{e}me 2.23, p. 144]{B-G-S2}.\smallskip }).\smallskip
\newline
(iii)} By \emph{(ii)} we may restrict ourselves to the case in which $%
X\left( N,\Delta _{\text{\emph{can}}}\right) $ has only canonical
singularities of index $1$ and, in particular, to cones $\tau \in \Delta _{%
\text{\emph{can}}}\left( 3\right) ,$ for which \emph{orb}$\left( \tau
\right) \in $ \emph{Sing}$\left( X\left( N,\Delta _{\text{\emph{can}}%
}\right) \right) .$ Each of these $\tau $'s has \emph{(}up to lattice
automorphism\emph{)} the form
\begin{equation*}
\tau _{P}=\left\{ \left( \lambda x_{1},\lambda x_{2},\lambda \right) \
\left| \ \lambda \in \mathbb{R}_{\geq 0}\right. ,\ \ \left(
x_{1},x_{2}\right) \in P\right\}
\end{equation*}
\emph{(}w.r.t. the standard rectangular lattice $\mathbb{Z}^{3}$\emph{)} for
some lattice polygon $P.$ There exists a composite of torus-equivariant,
crepant, projective, partial desingularizations
\begin{equation}
U_{\tau _{P}}\overset{f_{1}}{\longleftarrow }X_{\mathcal{S}_{1}}\overset{%
f_{2}}{\longleftarrow }X_{\mathcal{S}_{2}}\overset{f_{3}}{\longleftarrow }X_{%
\mathcal{S}_{3}}\longleftarrow \cdots \overset{f_{\kappa -1}}{\longleftarrow
}X_{\mathcal{S}_{\kappa -1}}\overset{f_{\kappa }}{\longleftarrow }X_{%
\mathcal{S}_{\kappa }}=Y_{\tau _{P}}  \label{FIRST SERIES}
\end{equation}
where $\mathcal{S}_{1},...,\mathcal{S}_{\kappa }$ are successive polygonal
subdivisions of $P$ \emph{(}with $X_{\mathcal{S}_{i}}$ the toric variety
associated to the cone supporting $\mathcal{S}_{i}$\emph{)}, such that each $%
f_{i}$ in \emph{(\ref{FIRST SERIES}) }is the usual toric blow-up of the
arising singular point and $Y_{\tau _{P}}$ is \emph{unique} w.r.t. this
property, \ possessing at most compound Du Val singularities whose types can
be written in a very short list\emph{\footnote{%
This list is the following:{\scriptsize
\begin{equation*}
\ \
\begin{tabular}{|c|c|c|}
\hline
\textbf{Cases} & $
\begin{array}{c}
\, \\
\text{\textbf{\ Possible cDu Val }} \\
\text{\textbf{singularities}} \\
\,
\end{array}
$ & $\text{\textbf{Types}}$ \\ \hline\hline
(i) & $
\begin{array}{c}
\  \\
\text{Spec}\left( \mathbb{C}[z_{1},z_{2},z_{3}]\,/\,\left(
z_{1}z_{2}-z_{3}^{\kappa }\right) \right) \times \text{Spec}\left( \mathbb{C}%
[z_{4}]\right) ,\ \ \kappa \in \mathbb{Z}_{\geq 2} \\
\
\end{array}
$ & $\mathbf{A}_{\kappa -1}\times \mathbb{C}$ \\ \hline
(ii) & $
\begin{array}{c}
\  \\
\text{Spec}\left( \mathbb{C}[z_{1},z_{2},z_{3},z_{4}]\,/\,\left(
z_{1}z_{2}-z_{3}^{\kappa }z_{4}^{\kappa +\lambda }\right) \right) ,\ \
\kappa \in \mathbb{Z}_{\geq 1},\,\lambda \in \mathbb{Z}_{\geq 0} \\
\
\end{array}
$ & $\mathbf{cA}_{2\kappa +\lambda -1}$ \\ \hline
(iii) & $
\begin{array}{c}
\  \\
\text{Spec}\left( \mathbb{C}[z_{1},z_{2},z_{3},z_{4}]\,/\,\left(
z_{1}z_{2}z_{3}-z_{4}^{2}\right) \right) \cong \text{Spec(}\mathbb{C}%
[t_{1},t_{2},t_{3}]^{G}) \\
\  \\
\ \text{(where }G\text{ is obtained by the linear representation} \\
\text{ of the Kleinian four-group into SL}(3,\mathbb{C}\text{))} \\
\
\end{array}
$ & $\mathbf{cD}_{4}$ \\ \hline
\end{tabular}
\end{equation*}
}}}. \emph{(}This is in fact an intrinsic procedure involving in each step
the auxiliary 2-dimensional subcones and the finite continued fraction
expansions for each vertex of the new central lattice polygon.\emph{%
)\smallskip }\ \newline
\emph{(iv) }There exists a composite of torus-equivariant, crepant,
projective, partial desingularizations
\begin{equation}
Y_{\tau _{P}}\overset{g_{1}}{\longleftarrow }X_{\widetilde{\mathcal{S}}_{1}}%
\overset{g_{2}}{\longleftarrow }X_{\widetilde{\mathcal{S}}_{2}}\overset{g_{3}%
}{\longleftarrow }X_{\widetilde{\mathcal{S}}_{3}}\longleftarrow \cdots
\overset{g_{\nu -1}}{\longleftarrow }X_{\widetilde{\mathcal{S}}_{\nu -1}}%
\overset{g_{\nu }}{\longleftarrow }X_{\widetilde{\mathcal{S}}_{\nu
}}=Z_{\tau _{P}}  \label{SECOND SERIES}
\end{equation}
of $Y_{\tau _{P}\text{,}}$ such that each $g_{j}$ in \emph{(\ref{SECOND
SERIES}) }is the usual toric blow-up of \ the ideal of the $1$-dimensional
locus of $X_{\widetilde{\mathcal{S}}_{j-1}}$ and $Z_{\tau _{P}}$ is \emph{%
unique} w.r.t. this property possessing at most \emph{ordinary double points}%
. \emph{(}$Z_{\tau _{P}}$ is a terminalization of \ $Y_{\tau _{P}}$\emph{)}.
\medskip \newline
\emph{(v)} All the $2^{\#\left( \text{\emph{ordinary double points}}\right)
} $ remaining choices to obtain ``full'' crepant desingularizations \emph{(}%
just by filling up our lattice triangulations by box diagonals\emph{)} are
realized by projective birational morphisms.
\end{theorem}

\noindent {}The proof involves explicit (constructive) techniques from
discrete geometry of cones, together with certain results of combinatorial
nature from \cite{FT-M}, \cite{Ishida}, \cite{ISHIDA-IWASHITA} and \cite
{Koelman}. (The projectivity of birational morphisms in (v) can be checked
by means of the ampelness criterion \ref{AMPLE}). Moreover, in \cite{DHZ},
it will be shown how toric geometry enables us to keep under control (in
each step) the Fano surfaces arising as exceptional prime divisors.

\begin{example}
As a simple example of how Theorem \ref{MAIN} (iii)-(v) works, take the cone
\begin{equation*}
\sigma =\text{pos}(P)\subset \mathbb{R}^{2}\times \left\{ 1\right\} \subset
\mathbb{R}^{3}
\end{equation*}
supporting the lattice tringle $P=$ conv$\left\{ \left( -3,3,1\right)
,\left( 3,1,1\right) ,\left( 0,-3,1\right) \right\} $ of Fig. \textbf{2}
(w.r.t. the standard rectangular lattice $\mathbb{Z}^{3}$).
\begin{figure}[h]
\begin{center}
\includegraphics[width=2in,height=2in]{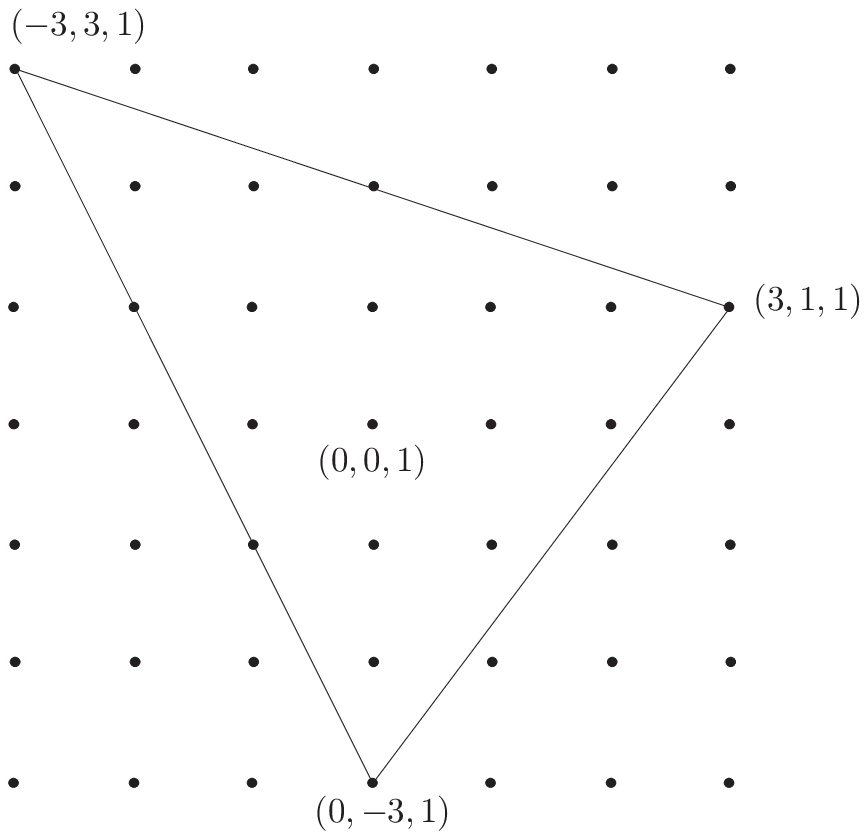} \vspace{0.3cm} \\[0pt]
\textbf{Fig. 2}
\end{center}
\end{figure}

\noindent {}\noindent The singularity $\left( U_{\sigma },\text{orb}\left(
\sigma \right) \right) $ is Gorenstein, and Proposition \ref{EMB} gives%
\footnote{%
This can be computed directly. A general formula, expressing the embedding
dimension in terms of the vertex coordinates of arbitrary initial lattice
polygons, will appear in \cite{DHZ}.}:
\begin{equation*}
\text{edim}\left( U_{\sigma },\text{orb}\left( \sigma \right) \right) =\#(%
\mathbf{Hilb}_{(\mathbb{Z}^{3})^{\vee }}\left( \sigma ^{\vee }\right) )=14.
\end{equation*}
Hence, by Theorem \ref{LAUFER-REID} (iv) and by Proposition \ref{PROP-REID}
(iii), the Laufer-Reid invariant of a general hyperplane section $\mathbb{H}$
through orb$\left( \sigma \right) $ equals
\begin{equation*}
\text{LRI}\left( \mathbb{H},\text{orb}\left( \sigma \right) \right) =\text{%
edim}\left( \mathbb{H},x\right) =\text{edim}\left( U_{\sigma },\text{orb}%
\left( \sigma \right) \right) -1=13.
\end{equation*}
Figure \textbf{3} shows the result of blowing up orb$\left( \sigma \right)
\in U_{\sigma }$ (equipped with the reduced subscheme structure).

\begin{figure}[h]
\begin{center}
\includegraphics[width=2.2in,height=2.2in]{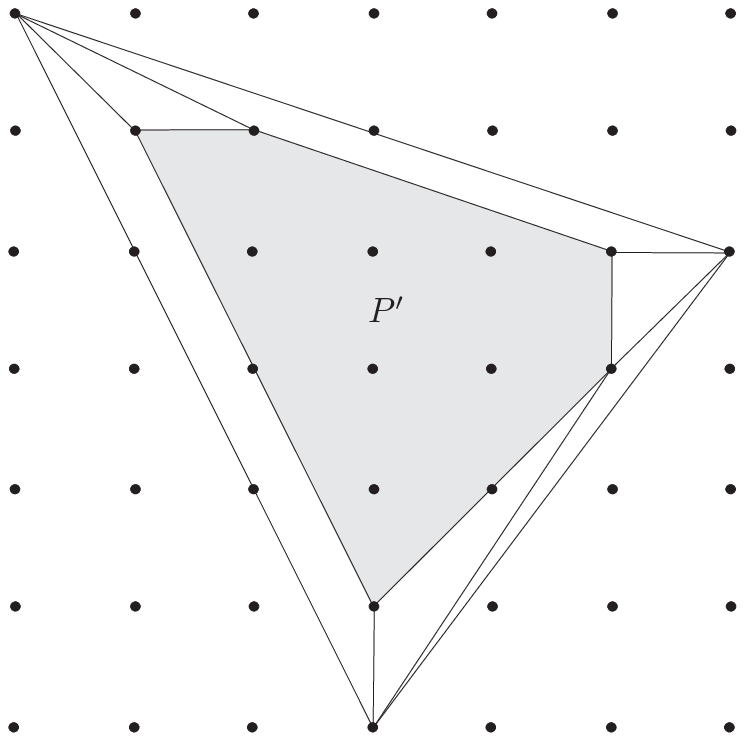} \vspace{0.3cm} \\[0pt]
\textbf{Fig. 3}
\end{center}
\end{figure}
\noindent {}It is worth mentioning that the ``central'' new subcone $\sigma
^{\prime }=$ pos$(P^{\prime })$ of $\sigma $ supports the lattice pentagon $%
P^{\prime }=$ conv$\left( \left\{ \left( -2,2,1\right) ,\left( -1,2,1\right)
,\left( 2,1,1\right) ,\left( 2,0,1\right) ,\left( 0,-2,1\right) \right\}
\right) $ which is nothing but the polygon defined as the convex hull of the
inner points of $P.$ Next, we perform \ref{MAIN} (iii)-(iv) again and again
until we arrive at the lattice polygonal subdivision of $P$ shown in Figure
\textbf{4}.
\begin{figure}[h]
\begin{center}
\includegraphics[width=2.2in,height=2.2in]{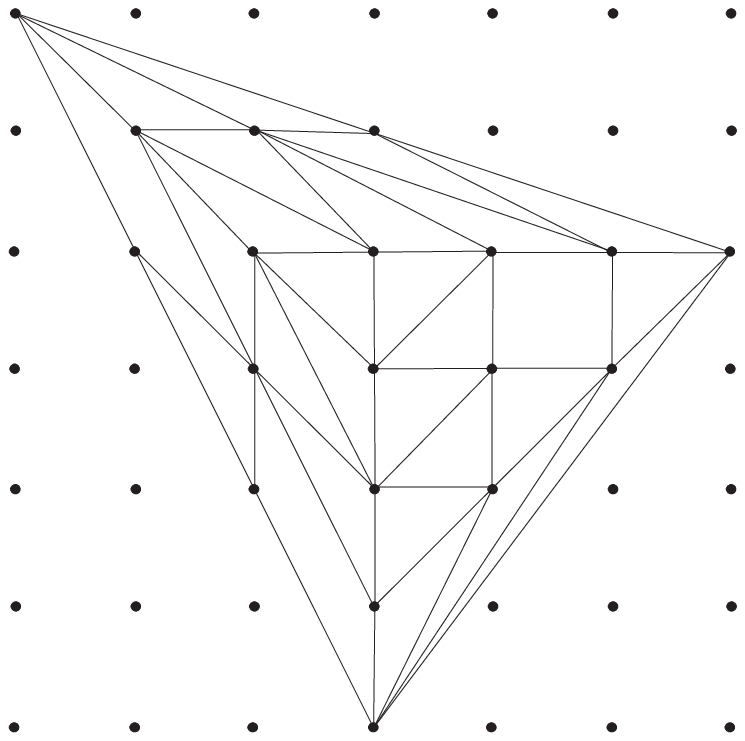} \vspace{0.3cm} \\[0pt]
\textbf{Fig. 4}
\end{center}
\end{figure}

\noindent \noindent {}In the last step (\ref{MAIN} (v)), all $2^{3}=8$
possible choices of completing the polygonal subdivision in Fig. \textbf{4 }%
to a triangulation (by filling up box diagonals) lead equally well to
crepant, projective, full resolutions of $\left( U_{\sigma },\text{orb}%
\left( \sigma \right) \right) .$
\end{example}

\begin{remark}
A method of how one may achieve ``$\mathbb{Q}$-factorialization'' of toric
singularities (after M. Reid \cite{REID3} and S. Mori) \textit{in arbitrary
dimensions} was partially discussed in Wi\'{s}niewski's lectures \cite
{WISNIEWSKI}.\medskip
\end{remark}

\noindent {}\textbf{Aknowledgements.} This paper is an expanded version of
my lecture notes which were distributed during the second week of the Summer
School ``Geometry of Toric Varieties'' in Grenoble (Institut Fourier, June
19 - July 7, 2000). In the first two talks, based on the preliminary notes
\cite{DAIS}, \cite{ZIEGLER}, I had the opportunity to give a sketch of proof
of Theorem \ref{MAIN} and to explain how one applies the corresponding
combinatorial mechanism in some concrete examples. (Complete proofs and
details will appear in the forthcoming joint work \cite{DHZ} with M. Henk
and G. M. Ziegler). I would like to express my warmest thanks to the
organizers M. Brion and L. Bonavero for the exemplary hospitality, to D. Cox
for useful remarks on much of the manuscript, as well as to several other
participants of the Summer School for many stimulating discussions.

\end{document}